\documentclass[11pt,twoside]{article}
\usepackage[plainpages=false]{hyperref}
\usepackage{amsfonts,latexsym,rawfonts,amsmath,amssymb,amsthm}
\usepackage{amsmath,amssymb,amsfonts,latexsym,lscape,rawfonts, color}
\numberwithin{equation}{section}

\newcommand{\pd}[2]{\frac {\partial #1}{\partial #2}}

\newcommand{\al}{\alpha}
\newcommand{\bb}{\beta}
\newcommand{\la}{\lambda}
\newcommand{\La}{\Lambda}
\newcommand{\oo}{\omega}

\newcommand{\dd}{\delta}
\newcommand{\Na}{\nabla}

\def\ga{\gamma}

\newcommand{\ee}{\epsilon}

\newcommand{\beq}{\begin{equation}}
\newcommand{\eeq}{\end{equation}}
\newcommand{\beqs}{\begin{eqnarray*}}
\newcommand{\eeqs}{\end{eqnarray*}}
\newcommand{\beqn}{\begin{eqnarray}}
\newcommand{\eeqn}{\end{eqnarray}}
\newcommand{\beqa}{\begin{array}}
\newcommand{\eeqa}{\end{array}}

\def\lra{\longrightarrow}
\def\La{\Lambda}
\def\la{\lambda}

\def\p{\prime}

\def\bc{\begin{center}}
\def\ec{\end{center}}

\def\p{\partial}
\def\b{\bar}

\def\cP{{\cal P}}

\def\cA{{\mathcal A}}

\def\cP{{\mathcal P}}

\def\RR{{\mathbb R}}
\def\NN{{\mathbb N}}

\def\CC{{\mathbb C}}

\def\begeq{\begin{equation}}
\def\endeq{\end{equation}}
\def\and{\quad{\rm and}\quad}

\def\Del{\Delta}

\let\lra=\longrightarrow

\def\mapright\#1{\,\smash{\mathop{\lra}\limits^{\#1}}\,}

\def\Vol{{\rm Vol}}

\def\ri{\rightarrow}
\def\pbp{\sqrt{-1}\partial\bar\partial}
\def\un{\underline}
\def\an{\; \;\;{ and}\;\;\;}

\newtheorem{prop}{Proposition}[section]
\newtheorem{theo}[prop]{Theorem}
\newtheorem{lem}[prop]{Lemma}
\newtheorem{claim}[prop]{Claim}

\newtheorem{rem}[prop]{Remark}

\newtheorem{defi}[prop]{Definition}

\title{{\bf\Large{On the K\"ahler-Ricci flow with small initial $E_1$ energy (I)}}}

\author{Xiuxiong Chen\footnote{The first named author is partially
supported by a NSF grant, while the third author was partially
supported by a NSF supplement grant.},  Haozhao Li and Bing Wang}

\date{}

\begin{document}
\bibliographystyle{plain}
\maketitle

\tableofcontents
\section{Introduction and main results }
\subsection{The motivation} In \cite{[chen-tian2]}
\cite{[chen-tian1]}, a family of functionals $E_k (k = 1,2,\cdots
n)$ was introduced by the first named author and G. Tian  to prove
the convergence of the K\"ahler-Ricci flow under appropriate
curvature assumptions. The aim of this program (cf. \cite{[chen1]})
is to study  how the lower bound of $E_1$ is used to derive the
convergence of the K\"ahler-Ricci flow, i.e., the existence of
K\"ahler-Einstein metrics.  We will address this question in
Subsection 1.2.   The corresponding problem of the relation between
the lower bound of $E_0$, which is the $K$-energy introduced by T.
Mabuchi,  and the existence of K\"ahler-Einstein metrics has been
extensively studied (cf. \cite{[BaMa]}, \cite{[chen-tian2]},
\cite{[chen-tian1]}, \cite{[Do]}). One interesting question in this
program is how the lower bound of $E_1$ compares to the lower bound
of $E_0$. We will give a satisfactory answer to this question
in Subsection 1.3.  \\

\subsection{The lower bound of $E_1$ and K\"ahler-Einstein metrics}
Let $(M, [\omega])$ be a polarized compact K\"ahler manifold with
$[\omega]=2\pi c_1(M)>0 $ (the first Chern class) in this paper. In
\cite{[chen1]}, the first named author proved a stability theorem of
the K\"ahler-Ricci flow near the infimum of $E_1$ under the
assumption that the initial metric has $Ric>-1$ and $|Rm|$ bounded.
Unfortunately, this stability theorem needs a topological assumption
\beq(-1)^n([c_1(M)]^{[2]}[\oo]^{[n-2]}-\frac
{2(n+1)}{n}[c_2(M)][\oo]^{[n-2]})\geq 0.\label{c1}\eeq The only
known compact manifold which satisfies this condition is $\CC P^n$,
which restricts potential applications of this result. The main
purpose of this paper is to remove this assumption.

\begin{theo}\label{main2}
Suppose that $M$ is pre-stable, and $E_1$ is bounded from below in $
[\oo]$. For any $\delta, \La>0,$ there exists a small positive
constant $\ee(\delta, \La)>0$ such that for any metric $\oo_0$ in
the subspace $\cA(\delta, \La, \oo, \ee)$ of K\"ahler metrics
$$\{\oo_{\phi}=\oo+\pbp \phi\;|\;
Ric(\oo_{\phi})>-1+\delta, \;  |Rm|(\oo_{\phi})\leq \La,\;
E_1(\oo_{\phi})\leq \inf E_1+\ee \},$$where $E_1(\oo')=E_{1,
\oo}(\oo')$, the K\"ahler-Ricci flow will deform it exponentially
fast to a K\"ahler-Einstein metric in the limit.
\end{theo}

\begin{rem} The condition that $M$ is pre-stable (cf. Defi \ref
{prestable}), roughly means that the complex structure doesn't jump
in the limit (cf. \cite{[chen1]}, \cite{[PhSt]}). In G. Tian's
definition of $K$-stability, this condition appears to be one of
three necessary conditions for a complex structure to be $K$-stable
(cf. \cite{[Do]}, \cite{[Tian2]}).
\end{rem}

\begin{rem}This gives a sufficient condition for the existence of  K\"ahler-Einstein metrics.
More interestingly, by a theorem of G. Tian \cite{[Tian2]}, this
also gives a sufficient condition for an algebraic manifold being
weakly K-stable.  One tempting question is: does this condition
imply weak K-stability directly?
\end{rem}

\begin{rem}
If we call the result in \cite{[chen1]} a ``pre-baby" step in this
ambitious program, then Theorem 1.1 and 1.5 should be viewed as a
``baby step" in this program. We wish to remove the assumption on
the bound of the bisectional curvature.  More importantly (cf.
Theorem 1.8 below), we wish to replace the condition on the Ricci
curvature in both Theorem 1.1 and 1.5 by a condition on the scalar
curvature. Then our theorem  really becomes a ``small energy" lemma.
\end{rem}

If we remove the condition of ``pre-stable", then

\begin{theo}\label{main} Suppose that $(M, [\oo])$ has no nonzero
holomorphic vector fields and $E_1$ is bounded from below in
$[\oo].$ For any $\delta, B,  \La>0,$ there exists a small
positive constant $\ee(\delta, B, \La, \oo)>0$ such that for any
metric $\oo_0$ in the subspace $\cA(\delta, B, \La, \ee)$ of
K\"ahler metrics
$$\{\oo_{\phi}=\oo+\pbp \phi\;|\; Ric(\oo_{\phi})>-1+\delta,\; |\phi|
\leq B, \;|Rm|(\oo_{\phi})\leq \La, \; E_1(\oo_{\phi})\leq \inf
E_1+\ee \}$$ the K\"ahler-Ricci flow will deform it exponentially
fast to a K\"ahler-Einstein metric in the limit.
\end{theo}

\begin{rem}  In light of Theorem \ref{main4} below, we can replace the
condition on $E_1$ by a corresponding  condition on $E_0.\;$
\end{rem}

\subsection{The relations between energy functionals $E_k$ }
Song-Weinkove  \cite{[SoWe]} recently  proved   that $E_k$ have a
lower bound
  on the space of K\"ahler metrics with nonnegative Ricci curvature
for K\"ahler-Einstein manifolds. Moreover, they also showed that
modulo holomorphic vector fields, $E_1$
  is proper if and only if there exists a K\"ahler-Einstein metric.
Shortly afterwards, N. Pali
    \cite{[Pali]} gave a formula between $E_1$ and the $K$-energy
$E_0$, which says that the $E_1$ energy
    is always bigger than the $K$-energy.
      Tosatti  \cite{[Tosatti]} proved that under some curvature assumptions,
    the critical point of $E_k$ is a K\"ahler-Einstein metric.
Pali's theorem means that $E_1$ has a lower bound if the $K$-energy
has a lower bound.
 A natural question is if
the converse holds. To our own surprise, we proved the following
result.

\begin{theo}\label{main4}$E_1$ is bounded from below if and only if
the $K$-energy is bounded from below  in the class $[\oo].$
Moreover, we have\footnote{For simplicity of notation, we will often
drop the subscript $\phi$ and write $|\nabla f|^2$ for $|\nabla
f|_{\phi}^2 $. But in an integral, $|\nabla f|^2$ is with respect to
the metric of the volume form. }
$$\inf_{\oo'\in [\oo]} E_{1, \oo}(\oo')=2\inf_{\oo'\in [\oo]} E_{0,
\oo}(\oo')-\frac 1{nV}\int_M\;|\Na h_{\oo}|^2 \oo^n,$$ where
$h_{\oo}$ is the Ricci potential function with respect to $\oo$.
\end{theo}

A crucial observation leads to this theorem is

\begin{theo} Along the K\"ahler-Ricci flow, $E_1$ will
decrease after finite time.
\label{main3}
\end{theo}

Theorem  \ref{main4} and \ref{main3} of course lead more questions
than they answer to.  For instance, is the properness of $E_k$
equivalent to the properness of $E_l$ for $k\neq l$? More subtlely,
are the corresponding notions of semi-stability ultimately
equivalent to each other? Is there a preferred functional in this
family, or a better linear combination of these $E_k$ functionals?
The first named author genuinely believes that
this observation opens the door for more interesting questions.   \\

Another interesting question is the relation of $E_k$ with various
notions of stability defined by algebraic conditions.  Theorem 1.1
and 1.5 suggest an indirect link of these functionals $E_k$ and
stability.   According to A. Futaki \cite{[futaki04]},  these
functionals may directly link to the asymptotic  Chow semi-stability
(note that the right hand side of (1.2) in \cite{[futaki04]} is
precisely equal to $\frac{d E_k}{dt}$ if one takes $p=k+1$ and $\phi
= c_1^k$, cf. Theorem 2.4 below).  It is highly interesting to
explore further in this direction.

\subsection{Main ideas of proofs of Theorem 1.1 and
1.5}

In \cite{[chen1]}, a topological condition is used to control the
$L^2$ norm of the bisectional curvature once the Ricci curvature is
controlled.   Using the parabolic Moser iteration arguments, this
gives a uniform bound on the full bisectional curvature.  In the
present paper, we need to find a new way to control the full
bisectional curvature under the flow.  The whole scheme of obtaining
this uniform estimate on curvatures depends on two crucial steps and
their dynamical interplay.\\

{\bf STEP 1:} The first step is to follow the approach of the
celebrated work of Yau on the Calabi conjecture (cf. {\cite{[Cao]}},
\cite{[Yau]}). The key point here is to control the $C^0$ norm of
the evolving K\"ahler potential $\phi(t)$, in particular, the growth
of $ u = {{\p\phi}\over {\p t}}$ along the K\"ahler-Ricci flow. Note
that $u$ satisfies
\[
{ {\p u}\over {\p t}} = \triangle_{\phi} u + u.
\]
Therefore, the crucial step is to control the first eigenvalue of
the Laplacian operator (assuming the traceless Ricci form is
controlled via an iteration process which we will describe as {\bf
STEP} 2 below).
  For our purpose, we need
to show that the first eigenvalue of the evolving Laplacian operator
is bigger than $1 + \ga$ for some fixed $\ga > 0.\;$ Such a problem
already appeared in \cite{[chen0]} and \cite{[chen-tian2]} since the
first eigenvalue of the Laplacian operator of K\"ahler-Einstein
metrics is exactly $1.\;$ If $Aut_r(M, J)\neq 0,\;$  the uniqueness
of K\"ahler-Einstein metrics implies that the dimension of the first
eigenspace is fixed; while the vanishing of the Futaki invariant
implies that $u(t)$ is essentially perpendicular to the first
eigenspace of the evolving metrics. These are two crucial
ingredients which allow one to squeeze out a small gap $\ga$ on the
first eigenvalue estimate. Following the approach in \cite {[chen0]}
and \cite{[chen-tian2]}, we can show that $u$ decays exponentially.
This in turn implies the $C^0$ bound on the evolving K\"ahler
potential. Consequently, this leads to control of all derivatives of
the evolving potential, in particular, the bisectional curvature. In
summary, as long as we have control of the first eigenvalue, one
controls the full bisectional curvature of the evolving K\"ahler
metrics.

For Theorem 1.5,  a crucial technique step is to use an estimate
obtained in \cite{[chenhe]} on the Ricci curvature tensor.\\

{\bf STEP 2:} Here we follow the Moser iteration techniques which
appeared in \cite{[chen1]}.  Assuming that the full bisectional
curvature is bounded by some large but fixed number, the norm of the
traceless bisectional curvature and the traceless Ricci tensor both
satisfy the following inequality:
\[
  {{\p u}\over {\p t}} \leq \triangle_{\phi} u + |Rm| u \leq
  \triangle_{\phi}
u + C\cdot u.
\]
If the curvature of the evolving metric is controlled in $L^p(p>
n)$, then the smallness of the energy $E_1$ allows you to control
the norm of the traceless Ricci tensor of the evolving K\"ahler
metrics (cf. formula (\ref{2})). According to Theorem \ref{lem2.8}
and \ref{theo4.18}, this will in turn give an improved estimate on
the first eigenvalue in a slightly longer period, but perhaps
without full uniform control of the bisectional curvature in the
``extra" time. However, this gives uniform control on the K\"ahler
potential which in turn gives the bisectional curvature in the
extended ``extra" time. We use the Moser iteration again to obtain
sharper control on the norm of the traceless bisectional
curvature.\\

Hence, the combination of the parabolic Moser iteration together
with Yau's estimate, yields the desired global estimate. In
comparison, the iteration process in \cite{[chen1]} is instructive
 and more direct.  The first named author believes that this
approach is perhaps more important than the mild results we obtained
there.

\subsection{The organization}  This paper is roughly organized as
follows: In Section 2, we review some basic facts in K\"ahler
geometry and necessary information on the K\"ahler-Ricci flow. We
also include some basic facts on the energy functionals $E_k$. In
Section 3, we prove Theorem \ref{main4} and \ref{main3}. In Section
4, we prove several technical theorems on the K\"ahler-Ricci flow.
The key results are the estimates of the first eigenvalue of the
evolving Laplacian, which are proved in Section 4.3. In Section 5, 6
we prove Theorem \ref{main2} and Theorem \ref{main}.\\


\noindent {\bf Acknowledgements}:  Part of this work was done while
the second named author was visiting University of Wisconsin-Madison
and he would like to express thanks for the hospitality. The second
named author would also like to thank
 Professor W. Y. Ding and X. H. Zhu for
their help and encouragement. The  first named author would like to
thank Professor P. Li of his interest and encouragement in this
project. The authors would like to thank the referees for numerous
suggestions which helped to improve the presentation. \vskip 1cm
\section{Setup and known results}
\subsection{Setup of notations}
Let $M$ be an $n$-dimensional compact K\"ahler manifold. A K\"ahler
metric can be given by its K\"ahler form $\omega$ on $M$. In local
coordinates $z_1, \cdots, z_n $, this $\omega$ is of the form
\[
\omega = \sqrt{-1} \displaystyle \sum_{i,j=1}^n\;g_{i
\overline{j}} d\,z^i\wedge d\,z^{\overline{j}}  > 0,
\]
where $\{g_{i\overline {j}}\}$ is a positive definite Hermitian
matrix function. The K\"ahler condition requires that $\omega$ is
a closed positive (1,1)-form. In other words, the following holds
\[
{{\partial g_{i \overline{k}}} \over {\partial z^{j}}} =
{{\partial g_{j \overline{k}}} \over {\partial z^{i}}}\qquad {\rm
and}\qquad {{\partial g_{k \overline{i}}} \over {\partial
z^{\overline{j}}}} = {{\partial g_{k \overline{j}}} \over
{\partial z^{\overline{i}}}}\qquad\forall\;i,j,k=1,2,\cdots, n.
\]
The K\"ahler metric corresponding to $\omega$ is given by
\[
\sqrt{-1} \;\displaystyle \sum_1^n \; {g}_{\alpha
\overline{\beta}} \; d\,z^{\alpha}\;\otimes d\, z^{
\overline{\beta}}.
\]
For simplicity, in the following, we will often denote by $\omega$
the corresponding K\"ahler metric. The K\"ahler class of $\omega$
is its cohomology class $[\omega]$ in $H^2(M,\RR).\;$ By the Hodge
theorem, any other K\"ahler metric in the same K\"ahler class is
of the form
\[
\omega_{\phi} = \omega + \sqrt{-1} \displaystyle \sum_{i,j=1}^n\;
{{\partial^2 \phi}\over {\partial z^i \partial z^{\overline{j}}}}
\;dz^i\wedge dz^{\bar j}
 > 0
\]
for some real valued function $\phi$ on $M.\;$ The functional space
in which we are interested (often referred to as the space of
K\"ahler potentials) is
\[
{\cal P}(M,\omega) = \{ \phi\in C^{\infty}(M, \RR) \;\mid\;
\omega_{\phi} = \omega + \sqrt{-1} {\partial} \overline{\partial}
\phi > 0\;\;{\rm on}\; M\}.
\]
Given a K\"ahler metric $\omega$, its volume form  is
\[
\omega^n = n!\;\left(\sqrt{-1} \right)^n \det\left(g_{i
\overline{j}}\right) d\,z^1 \wedge d\,z^{\overline{1}}\wedge \cdots
\wedge d\,z^n \wedge d \,z^{\overline{n}}.
\]
Its Christoffel symbols are given by
\[
\Gamma^k_{i\,j} = \displaystyle \sum_{l=1}^n\;g^{k\overline{l}}
{{\partial g_{i \overline{l}}} \over {\partial z^{j}}} ~~~{\rm
and}~~~ \Gamma^{\overline{k}}_{\overline{i} \,\overline{j}} =
\displaystyle \sum_{l=1}^n\;g^{\overline{k}l} {{\partial g_{l
\overline{i}}} \over {\partial z^{\overline{j}}}},
\qquad\forall\;i,j,k=1,2,\cdots n.
\]
The curvature tensor is
\[
R_{i \overline{j} k \overline{l}} = - {{\partial^2 g_{i \overline
{j}}} \over {\partial z^{k} \partial z^{\overline{l}}}} +
\displaystyle \sum_ {p,q=1}^n g^{p\overline{q}} {{\partial g_{i
\overline{q}}} \over {\partial z^{k}}}  {{\partial g_{p
\overline{j}}} \over {\partial z^{\overline{l}}}},
\qquad\forall\;i,j,k,l=1,2,\cdots n.
\]
We say that $\omega$ is of nonnegative bisectional curvature if
\[
R_{i \overline{j} k \overline{l}} v^i v^{\overline{j}} w^k w^
{\overline{l}}\geq 0
\]
for all non-zero vectors $v$ and $w$ in the holomorphic tangent
bundle of $M$. The bisectional curvature and the curvature tensor
can be mutually determined. The Ricci curvature of $\omega$ is
locally given by
\[
R_{i \overline{j}} = - {{\partial}^2 \log \det (g_{k
\overline{l}}) \over {\partial z_i \partial \bar z_j }} .
\]
So its Ricci curvature form is
\[
{\rm Ric}({\omega}) = \sqrt{-1} \displaystyle \sum_{i,j=1}^n \;R_{i
\overline{j}}\; d\,z^i\wedge d\,z^{\overline{j}} = -\sqrt{-1}
\partial \overline {\partial} \log \;\det (g_{k \overline{l}}).
\]
It is a real, closed (1,1)-form. Recall that $[\omega]$ is called a
canonical K\"ahler class if this Ricci form is cohomologous to
$\lambda \,\omega$ for some constant $\lambda$\,. In our setting, we
require $\lambda = 1.\;$

\subsection{The K\"ahler-Ricci flow}

Let us assume that the first Chern class  $c_1(M)$ is positive.
Choose an initial K\"ahler metric $\omega$ in $2\pi c_1(M).\;$ The
normalized K\"ahler-Ricci flow (cf. \cite{[Ha82]}) on a K\"ahler
manifold $M$ is of the form
\begin{equation}
{{\partial g_{i \overline{j}}} \over {\partial t }} = g_{i
\overline{j}} - R_{i \overline{j}}, \qquad\forall\; i,\; j=
1,2,\cdots ,n. \label{eq:kahlerricciflow}
\end{equation}
It follows that on the level of K\"ahler potentials, the Ricci
flow becomes
\begin{equation}
{{\partial \phi} \over {\partial t }} =  \log {{\omega_{\phi}}^n
\over {\omega}^n } + \phi - h_{\omega} , \label{eq:flowpotential}
\end{equation}
where $h_{\omega}$ is defined by
\[
{\rm Ric}({\omega})- \omega = \sqrt{-1} \partial \overline{\partial}
h_ {\omega}, \; {\rm and}\;\displaystyle \int_M\; (e^{h_{\omega}} -
1)  {\omega}^n = 0.
\]
Then the evolution equation for bisectional curvature is

\begin{eqnarray}{{\partial }\over {\partial t}} R_{i \overline{j} k
\overline{l}} & = & \bigtriangleup R_{i \overline{j} k
\overline{l}} + R_{i \overline{j} p \overline{q}} R_{q
\overline{p} k \overline{l}} - R_{i \overline{p} k \overline{q}}
R_{p \overline{j} q \overline{l}} + R_{i \overline{l} p
\overline{q}} R_{q \overline{p} k \overline{j}} + R_{i
\overline{j} k \overline{l}} \nonumber\\
& &  \;\;\; -{1\over 2} \left( R_{i \overline{p}}R_{p \overline{j}
k \overline{l}}  + R_{p \overline{j}}R_{i \overline{p} k
\overline{l}} + R_{k \overline{p}}R_{i \overline{j} p
\overline{l}} + R_{p \overline{l}}R_{i \overline{j} k
\overline{p}} \right). \label{eq:evolutio of curvature1}
\end{eqnarray}
Here $\Delta$ is the Laplacian of the metric $g(t).$ The evolution
equation for Ricci curvature and scalar curvature are
\begin{eqnarray} {{\p R_{i \b j}}\over {\p t}} & = & \triangle
R_{i\b j} + R_{i\b j p \b q} R_{q \b p} -R_{i\b p} R_{p \b
j},\label
{eq:evolutio of curvature2}\\
{{\p R}\over {\p t}} & = & \triangle R + R_{i\b j} R_{j\b i}- R.
\label{eq:evolutio of curvature3}
\end{eqnarray}

For direct computations and using the evolving frames, we can obtain
the following evolution equations for the bisectional curvature:
\begin{equation}
\pd {R_{i\bar jk\bar l}}{t} =\Delta R_{i\bar jk\bar l}- R_{i\bar
jk\bar l}+R_{i \bar j m\bar n}R_{n\bar m k\bar l}-R_{i\bar m k\bar
n}R_{m\bar j n\bar l}+R_{i\bar l m\bar n}R_{n\bar m k\bar l}
\label{eq:evolution of curvature4}
\end{equation}

As usual, the flow equation (\ref{eq:kahlerricciflow}) or
(\ref{eq:flowpotential}) is referred to as the K\"ahler-Ricci flow
on $M$. It was proved by Cao \cite{[Cao]}, who followed Yau's
celebrated work \cite{[Yau]}, that the K\"ahler-Ricci flow exists
globally for any smooth initial K\"ahler metric. It was proved by S.
Bando \cite{[Bando]} in dimension 3 and N. Mok \cite{[Mok]} in all
dimensions that the positivity of the bisectional curvature is
preserved under the flow. In \cite{[chen-tian2]} and
\cite{[chen-tian1]}, the first named author and G. Tian proved that
the K\"ahler-Ricci flow, in a K\"ahler-Einstein manifold, initiated
from a metric with positive bisectional curvature converges to a
K\"ahler-Einstein metric with constant bisectional curvature.  In
unpublished work on the K\"ahler-Ricci flow, G. Perelman proved,
along with other results, that the scalar curvature is always
uniformly bounded.

\subsection{Energy functionals $E_k$}
In \cite{[chen-tian2]}, a family of energy functionals $E_k (k=0, 1,
2,\cdots, n)$ was introduced and these functionals played an
important role there. First, we recall the definitions of these
functionals.
\begin{defi} For any $k=0, 1, \cdots, n,$ we define a functional
$E_k^0$ on $\cP(M, \oo)$ by
$$E_{k, \oo}^0(\phi)=\frac 1V\int_M\; \Big(\log \frac {\oo^n_{\phi}}
{\oo^n}-h_{\oo}\Big)\Big( \sum_{i=0}^k\; Ric(\oo_{\phi})^i\wedge
\oo^{k-i}\Big)\wedge \oo_{\phi}^{n-k}+\frac 1V\int_M\; h_{\oo}\Big(
\sum_{i=0}^k\; Ric(\oo)^i\wedge \oo^{k-i}\Big)\wedge \oo^{n-k}.$$
\end{defi}

\begin{defi}For any $k=0, 1, \cdots, n$, we define $J_{k, \oo}$ as
follows
$$J_{k, \oo}(\phi)=-\frac {n-k}{V}\int_0^1\;\int_M\;\pd {\phi(t)}{t}(\oo_{\phi(t)}
^{k+1}-\oo^{k+1})\wedge \oo_{\phi(t)}^{n-k-1} \wedge dt,$$ where
$\phi(t)(t\in [0, 1])$ is a path from $0$ to $\phi$ in $\cP(M,
\oo)$.
\end{defi}

\begin{defi}For any $k=0, 1, \cdots, n,$ the functional $E_{k, \oo}$
is defined as follows
$$E_{k, \oo}(\phi)=E_{k, \oo}^0(\phi)-J_{k, \oo}(\phi).$$
For simplicity, we will often drop the subscript $\oo$.
\end{defi}

By direct computation, we have
\begin{theo} For any $k=0, 1, 2,\cdots, n,$ we have
\beqs \frac {dE_k}{dt}=&&\frac {k+1}{V}\int_M\;\Delta_{\phi}\dot\phi
Ric(\oo_{\phi})^k\wedge\oo_{\phi}^{n-k}\\&&-\frac
{n-k}{V}\int_M\;\dot\phi
(Ric(\oo_{\phi})^{k+1}-\oo_{\phi}^{k+1})\wedge \oo_{\phi}^{n-k-1}.
\eeqs Here $\phi(t)$ is any path in $\cP(M, \oo)$.
\end{theo}

\begin{rem}Note that
$$
\frac {dE_0}{dt}=-\frac nV\int_M\; \dot\phi
(Ric(\oo_{\phi})-\oo_{\phi})\wedge \oo_{\phi}^{n-1}.$$ Thus, $E_0$
is the well-known $K$-energy.
\end{rem}
\begin{theo} Along the K\"ahler-Ricci flow where $Ric(\oo_{\phi})>-\oo_{\phi}$ is
preserved, we have
$$\frac {dE_k}{dt}\leq-\frac {k+1}{V}\int_M\;(R(\oo_{\phi})-r) Ric(\oo_{\phi})^k\wedge\oo_{\phi}^{n-k}.$$
 When $k=0, 1,$ we have \beqs
\frac {dE_0}{dt}&=&-\frac 1V\int_M\;|\Na \dot\phi|^2\oo_{\phi}^n
\leq 0,\\
\frac {dE_1}{dt}&\leq&-\frac
2V\int_M\;(R(\oo_{\phi})-r)^2\oo_{\phi}^n\leq 0. \eeqs
\end{theo}

Recently, Pali in \cite{[Pali]} found the following formula, which
will be used in this paper.
\begin{theo}\label{pali}For any $\phi\in \cP(M, \oo)$, we have
$$E_{1, \oo}(\phi)=2E_{0, \oo}(\phi)+\frac 1{nV}\int_M\; |\Na u|^2
\oo_{\phi}^n -\frac 1{nV}\int_M\; |\Na h_{\oo}|^2 \oo^n,$$ where
$$u=\log\frac {\oo^n_{\phi}}{\oo^n}+\phi-h_{\oo}.$$
\end{theo}

\begin{rem}This formula directly implies that if $E_0$ is bounded
from below, then $E_1$ is bounded from below on $\cP(M, \oo)$.
\end{rem}

\begin{rem}In a forthcoming paper \cite{[chenli]}, the second named author will
generalize Theorem \ref{pali} to all the functionals $E_k(k\geq 1)$,
and discuss some interesting relations between $E_k.$
\end{rem}
\section{Energy functionals $E_0$ and $E_1$}
In this section, we want to prove Theorem \ref{main4} and
\ref{main3}.

\subsection{Energy functional $E_1$ along the K\"ahler-Ricci flow}
The following theorem  is well-known in literature  (cf.
\cite{[chen2]}).

\begin{lem}\label{sec2}The minimum of the scalar curvature along the K\"ahler-Ricci flow,
if negative, will increase to zero
exponentially.
\end{lem}
\begin{proof}  Let $\mu(t)=-\min_M R(x, 0)e^{-t},$ then
\beqs\pd {}t(R+\mu(t))&=&\Delta (R+\mu(t))+|Ric|^2-(R+\mu(t))\\
&\geq &\Delta (R+\mu(t))-(R+\mu(t)). \eeqs Since $R(x,
0)+\mu(0)\geq 0$, by the maximum principle we have
$$R(x, t)\geq -\mu(t)=\min_M R(x, 0)e^{-t}.$$
\end{proof}

Using the above lemma, the following theorem is an easy corollary
of Pali's formula.
\begin{theo}\label{thm3.2}Along the K\"ahler-Ricci flow, $E_1$ will
decrease after finite time. In particular, if the initial scalar
curvature $R(0)>-n+1$, then there is a small constant $\dd>0$
depending only on
  $n$ and $\min_{x\in M} R(0)$ such that for all time $t>0$, we have
  \beq \frac {d}{dt}E_1\leq -\frac {\dd}V \int_M\;|\Na \dot \phi|^2
\oo_{\phi}^n\leq
  0.\label{aaa}\eeq

\end{theo}
\begin{proof} Along the K\"ahler-Ricci flow, the evolution equation
for $|\Na \dot \phi|^2$ is
$$\pd {}t|\Na \dot \phi|^2=\Delta_{\phi} |\Na \dot \phi|^2-|\Na \Na \dot
\phi|^2-|\Na \bar \Na\dot \phi|^2+|\Na \dot \phi|^2.$$ By Theorem
\ref{pali}, we have \beqs\frac {d}{dt}E_1&=&-\frac 2V\int_M\;|\Na
\dot \phi|^2\oo_{\phi}
^n +\frac 1{nV}\frac {d}{dt}\int_M\; |\Na \dot \phi|^2\oo_{\phi}^n\\
&= &-\frac 2V\int_M\;|\Na \dot \phi|^2\oo_{\phi}^n +\frac
1{nV}\int_M\;(-|\Na \Na \dot \phi|^2-|\Na \bar \Na\dot
\phi|^2+|\Na \dot \phi|^2+|\Na \dot \phi|^2(n-R))\oo_{\phi}^n
\eeqs If the scalar curvature at the initial time $R(x, 0)\geq
-n+1+n\dd (n\geq 2)$ for some small $\dd>0$,  by Lemma \ref{sec2}
for all time $t>0$ we have $R(x, t)\geq -n+1+n\dd.$ Then we have
\beqn\frac {d}{dt}E_1 &\leq&-\frac 2V\int_M\;|\Na \dot
\phi|^2\oo_{\phi}^n +\frac 1{nV}\int_M\;(-|\Na \Na \dot \phi|^2-
|\Na \bar \Na\dot \phi|^2+|\Na \dot \phi|^2+|\Na \dot
\phi|^2(2n-1-n
\delta))\oo_{\phi}^n\nonumber\\
&\leq&-\frac {\dd}V \int_M\;|\Na \dot \phi|^2\oo_{\phi}^n.
\label{aeq1}\eeqn Otherwise, by Lemma \ref{sec2} after finite time
$T_0=\log\frac {-\min_M R(x, 0)}{n-1-n\dd},$  we still have $R(x,
t)>-n+1+n\delta$ for small $\dd
 >0.$ Thus, the inequality (\ref{aeq1}) holds.

If $n=1$, by direct calculation we have
$$\frac {dE_1}{dt}=-\frac 2V\int_M\; (R(\oo_{\phi})-1)^2\oo_{\phi}=-\frac
2V\int_M\; (\Delta_{\phi}\dot\phi)^2\oo_{\phi}.$$ If the initial
scalar curvature $R(0)>0$, by R. Hamilton's results in \cite{[Ha88]}
the scalar curvature has a uniformly positive lower bound. Thus,
$R(t)\geq c>0$ for some constant $c>0$ and all $t>0.$ Therefore, by
the proof of Lemma \ref{lem2.13} in section \ref{section4.4.1} the
first eigenvalue of $\Delta_{\phi}$ satisfies $\la_1(t)\geq c.$ Then
$$\frac {dE_1}{dt}=-\frac 2V\int_M\; (\Delta_{\phi}\dot\phi)^2\oo_
{\phi}\leq -\frac {2c}V\int_M\; |\Na\dot\phi|^2\oo_{\phi}.$$ The
theorem is proved.

\end{proof}

\subsection{On the lower bound of $E_0$ and $E_1$} In this section, we
will prove Theorem \ref{main4}.  Recall the generalized energy:
\beqs I_{\oo}(\phi)&=&\frac 1V\sum_{i=0}^{n-1}\int_M\; \sqrt{-1}
\partial \phi\wedge\bar \partial \phi\wedge \oo^i\wedge\oo_{\phi}^
{n-1-i},\\
J_{\oo}(\phi)&=&\frac 1V\sum_{i=0}^{n-1}\frac {i+1}{n+1}\int_M\;
\sqrt{-1}\partial \phi\wedge\bar \partial \phi\wedge
\oo^i\wedge\oo_{\phi}^{n-1-i}. \eeqs By direct calculation, we can
prove
$$0\leq I-J\leq I\leq (n+1)(I-J)$$
and for any K\"ahler potential $\phi(t)$
$$\frac {d}{dt}(I_{\oo}-J_{\oo})(\phi(t))=-\frac
1V\int_M\;\phi\Delta_{\phi} \dot\phi \;\oo_{\phi}^n.$$

The behaviour of $E_1$ for the family of K\"ahler potentials
$\phi(t)$ satisfying the equation (\ref{eq1}) below has been studied
by Song and Weinkove in \cite{[SoWe]}. Following their ideas, we
have the following lemma.

\begin{lem}\label{lema3.3}For
any K\"ahler metric $\oo_0\in[\oo]$, there exists a K\"ahler
metric $\oo_0'\in [\oo]$ such that $Ric(\oo_0')>0$ and
$$E_0(\oo_0)\geq E_0(\oo_0').$$
\end{lem}

\begin{proof}We consider the complex Monge-Amp\`{e}re equation
\beq(\oo_0+\pbp \varphi)^n=e^{th_{0}+c_t}\oo_0^n,\label{eq1}\eeq
where $h_{0}$ satisfies the following equation
$$Ric(\oo_0)-\oo_0=\pbp h_{0},\qquad \frac 1V\int_M\;e^{h_0}
\oo_0^n=1$$ and $c_t$ is the constant chosen so that
$$\int_M\; e^{th_0+c_t}\oo_0^n=V.$$
By Yau's results in \cite{[Yau]}, there exists a unique
$\varphi(t)(t\in [0, 1])$ to the equation (\ref{eq1}) with $\int_M\;
\varphi \oo_0^n=0.$ Then $\varphi(0)=0.$ Note that the equation
(\ref{eq1}) implies \beq Ric(\oo_{\varphi})=\oo_{\varphi}+(1-t)\pbp
h_0-\pbp\varphi,\label{eq2}\eeq and
$$\Delta_{\varphi}\dot\varphi=h_0+c_t'.$$ By the definition of $E_0$
we have \beqs \frac {d}{dt}E_0(\varphi(t))&=&-\frac 1V\int_M\;
\dot\varphi(R(\oo_{\varphi})-n)\oo_{\varphi}^n\\
&=&-\frac 1V\int_M\;
\dot\varphi((1-t)\Delta_{\varphi}h_0-\Delta_{\varphi}\varphi)\oo_
{\varphi}^n\\
&=&-\frac {1-t}{V}\int_M\;\Delta_{\varphi} \dot\varphi h_0\;
\oo_{\varphi}^n+\frac
1V\int_M\;\varphi\Delta_{\varphi}\dot\varphi\oo_{\varphi}^n\\
&=&-\frac {1-t}{V}\int_M\;(\Delta_{\varphi} \dot\varphi)^2\;\oo_
{\varphi}^n-\frac {d}{dt}(I-J)_{\oo_0}(\varphi). \eeqs Integrating
the above formula from $0$ to $1$, we have
$$E_0(\varphi(1))-E_0(\oo_0)=-\frac 1V\int_0^1(1-s)\int_M\;(\Delta_
{\varphi} \dot\varphi)^2\; \oo_{\varphi}^n\wedge
ds-(I-J)_{\oo_0}(\varphi(1))\leq 0.$$ By the equation (\ref{eq2}),
we know $Ric(\oo_{\varphi(1)})> 0$. This proves the lemma.
\end{proof}

Now we can prove Theorem \ref{main4}.
\begin{theo}$E_1$ is bounded from below if and only if the $K$-energy
is bounded from below  in the class $[\oo].$ Moreover, we have
$$\inf_{\oo'\in [\oo]} E_{1}(\oo')=2\inf_{\oo'\in [\oo]} E_{0}(\oo')-
\frac 1{nV}\int_M\; |\Na h_{\oo}|^2\oo^n.$$
\end{theo}

\begin{proof}It is sufficient to show that if $E_1$ is bounded from below,
then $E_0$ is bounded from below. For any K\"ahler metric $\oo_0$,
by Lemma \ref{lema3.3} there exists a K\"ahler metric
$\oo_0'=\oo+\pbp \varphi_0$ such that
$$Ric(\oo_0')\geq c>0,\qquad E_0(\oo_0)\geq E_0(\oo_0'),$$
where $c$ is a constant depending only on $\oo_0.$  Let $\varphi(t)$
be the solution to the K\"ahler-Ricci flow with the initial metric
$\oo_0',$
$$\pd {\varphi}{t}=\log\frac {\oo_{\varphi}^n}{\oo^n}+\varphi-h_{\oo},
\qquad \varphi(0)=\varphi_0.$$ Then for any $t>s\geq0$, by Theorem
\ref{thm3.2} we have \beq E_1(t)-E_1(s)\leq
2\dd(E_0(t)-E_0(s)),\label{qeq3}\eeq where $E_1(t)=E_1(\oo,
\oo_{\varphi(t)})$ and $\dd=\frac {n-1}{2n}$ if $n\geq 2,$ or
$\dd=c>0 $ if $n=1$. Here $c$ is a constant obtained in the proof of
Theorem \ref{thm3.2}. By Theorem \ref{pali} we have
$$E_1(t)-2E_0(s)-\frac 1{nV}\int_M\;|\Na\dot\varphi|^2\oo_{\varphi}^n
(s)+C_{\oo}\leq \dd(E_1(t) -\frac
1{nV}\int_M\;|\Na\dot\varphi|^2\oo_{\varphi}^n(t)+C_{\oo})-2\dd
E_0(s).$$ i.e. \beq E_1(t)-\frac
1{n(1-\dd)V}\int_M\;|\Na\dot\varphi|^2\oo_{\varphi}^n(s)+ \frac
{\dd}{n(1-\dd)V}\int_M\;|\Na\dot\varphi|^2\oo_{\varphi}^n(t)+C_{\oo}\leq
2E_0(s),\label{qeq4}\eeq where $C_{\oo}=\frac 1{nV}\int_M\; |\Na
h_{\oo}|^2\oo^n.$ By (\ref{qeq3}) we know $E_0$ is bounded from
below along the K\"ahler-Ricci flow. Thus there exists a sequence of
times $t_m$ such that
$$\int_M\;|\Na\dot\varphi|^2\oo_{\varphi}^n(t_m)\ri 0,\qquad m\ri
\infty.$$ We choose $t=t_m$ and let $m\ri \infty$ in (\ref{qeq4}),
$$\inf E_1-\frac
1{n(1-\dd)V}\int_M\;|\Na\dot\varphi|^2\oo_{\varphi}^n(s)+C_{\oo}\leq
2E_0(s)\leq 2E_0(\oo_0'),
$$ where the last inequality is because $E_0$ is decreasing along
the K\"ahler-Ricci flow. Thus, we choose $s=t_m$ again and let $m\ri
\infty$,
$$\inf E_1+C_{\oo}\leq 2E_0(\oo_0')\leq 2E_0(\oo_0).$$
Thus, $E_0$ is bounded from below in $[\oo]$, and
$$\inf E_1\leq  2\inf E_0-C_{\oo}.$$
On the other hand, for any $\oo'\in [\oo]$ we have
$$E_1(\oo')\geq 2E_0(\oo')-C_{\oo}.$$
Combining the last two inequalities, we have $\inf E_1=2\inf
E_0-C_{\oo}$. Thus, the theorem is proved.
\end{proof}

\section{Some technical lemmas}
In this section, we will prove some technical lemmas, which will be
used in the proof of Theorem \ref{main2} and \ref{main}. These
lemmas are based on the K\"ahler-Ricci flow
$$\pd {\phi}{t}=\log\frac {\oo_{\phi}^n}{\oo^n}+\phi-h_{\oo}.$$
Most of these results are taken from
\cite{[chen1]}-\cite{[chen-tian1]}. The readers are referred to
these papers  for the details. Here we will prove some of them for
completeness.

\subsection{Estimates of the Ricci curvature}
The following result shows that we can control the curvature
tensor in a short time.

\begin{lem}\label{lem2.1}(cf. \cite{[chen1]})
Suppose that for some $\dd>0$, the curvature of $\oo_0=\oo+\pbp
\phi(0)$ satisfies the following conditions
$$\left \{\begin{array}{lll}|Rm|(0)&\leq& \La,\\
R_{i\bar j}(0)&\geq& -1+\dd.
\end{array}\right. $$
Then there exists a constant $T(\dd, \La)>0$, such that  for the
evolving K\"ahler metric $\oo_t(0\leq t\leq 6T)$, we have the
following
\beq \left \{\begin{array}{lll}|Rm|(t)&\leq& 2\La,\\
R_{i\bar j}(t)&\geq& -1+\frac {\dd}2.
\end{array}\right. \label{1}\eeq
\end{lem}

\begin{lem}

\label{lem2.2}(cf. \cite{[chen1]})If $E_1(0)\leq \inf_{\oo'\in
[\oo] }E_1(\oo')+\ee,$  and
$$Ric(t)+\oo(t)\geq \frac {\dd}2>0,\qquad \forall t\in
[0, T],$$ then along the K\"ahler-Ricci flow we have \beq\frac
1V\int_0^{T}\int_M\;\;|Ric-\oo|^2(t)\oo_ {\phi}^n\wedge dt\leq \frac
{\ee}2.\label{2}\eeq
\end{lem}

Since we have the estimate of the Ricci curvature, the following
theorem shows that the Sobolev constant is uniformly bounded if
$E_1$ is small.

\begin{prop}\label{lem2.3}(cf. \cite{[chen1]}) Along the K\"ahler-Ricci flow, if $E_1(0)\leq \inf_{\tilde\oo\in [\oo]
}E_1(\tilde\oo)+\ee,$ and for any $t\in [0, T],$
$$Ric(t)+\oo(t)\geq 0, $$ the diameter of the
evolving metric $\oo_{\phi}$ is uniformly bounded for $t\in [0, T].$
As $\ee\ri 0, $ we have $D\ri \pi.$ Let $\sigma(\ee)$ be the maximum
of the Sobolev and Poincar\'e constants with respect to the metric
$\oo_{\phi}$. As $\ee\ri 0,$ we have $\sigma(\ee)\leq
\sigma<+\infty.$ Here $\sigma$ is a constant independent of $\ee.$
\end{prop}

Next we  state a parabolic version of Moser iteration argument (cf.
\cite{[chen-tian1]}).

\begin{prop}\label{lem2.17} Suppose the Sobolev and Poincar\'e
constants of the evolving K\"ahler metrics $g(t)$ are both uniformly
bounded by $\sigma$. If a nonnegative function $u$ satisfies the
following inequality
$$\pd {}tu\leq \Del_{\phi} u+f(t, x)u, \;\; \forall \,t\in (a, b),$$
where $|f|_{L^p(M, g(t))}$ is uniformly bounded by some constant $c$
for some $p>\frac m2 $, where $m=2n=\dim_{\RR}M$, then for any
$\tau\in (0, b-a)$ and any $t\in (a+\tau, b)$, we have\footnote{The
constant $C$ may differ from line to line. The notation $C(A, B,
...)$ means that the constant $C$ depends only on $A, B, ...$.}
$$u(t)\leq \frac {C(n, \sigma, c)}{\tau^{\frac {m+2}{4}}}\Big(\int_{t-
\tau}^t\int_M \;u^2 \,\oo_{\phi}^n\wedge ds\Big)^{\frac 12}.$$
\end{prop}

By the above Moser iteration, we can show the following lemma.
\begin{lem}\label{lem2.5} For any $\dd, \La>0,$ there exists a small
positive constant $\ee(\dd, \La)>0$ such that if the initial metric
$\oo_0$ satisfies the following condition: \beq Ric(0)>-1+\dd,\;
|Rm(0)|\leq \La, \;E_1(0)\leq\inf E_1+\ee,\label{z1}\eeq then after
time $2T$ along the K\"ahler-Ricci flow, we have
\beq|Ric-\oo|(t)\leq C_1(T, \La)\ee, \qquad\forall t\in [2T,
6T]\label{2.5}\eeq and \beq|\dot\phi-c(t)|_{C^0}\leq C(\sigma)C_1(T,
\La)\ee, \qquad\forall t \in [2T, 6T],\label{z2}\eeq where $c(t)$ is
the average of $\dot\phi$ with respect to the metric $g(t)$, and
$\sigma$ is the uniformly upper bound of the Sobolev and Poincar\'e
constants in Proposition \ref{lem2.3}.
\end{lem}
\begin{proof} Let $Ric^0=Ric-\oo.$ Then $u=|Ric^0|^2(t)$ satisfies the
parabolic inequality
$$\pd ut\leq \Delta_{\phi} u+c(n)|Rm|_{g(t)}u,$$
Note that by Lemma \ref{lem2.1}, $|Rm|(t)\leq 2\La,$ for $0\leq
t\leq 6T.$  Then applying Lemma \ref{lem2.1} again and Lemma
\ref{lem2.2} for $t\in [2T, 6T]$, we have \footnote{Since the volume
$V$ of the K\"ahler manifold $M$ is fixed for the  metrics in the
same K\"ahler class, the constant $C(T, \La)$ below should depend on
$V$, but we don't specify this for simplicity.} \beqs
|Ric^0|^2(t)&\leq &C(\La, T)\Big(\int_{0}^{6T}\int_M\;\;|Ric-\oo|^4(t)\oo_{\phi}^n\wedge dt\Big)^{\frac 12}\\
&\leq &C(\La, T)(1+\La)\Big(\int_{0}^{6T}\int_M\;\;|Ric-\oo
|^2(t)\oo_{\phi}^n\wedge dt\Big)^{\frac
12}\\
&\leq &C(\La, T)\sqrt{\ee}. \eeqs Thus, \beq|Ric-\oo|(t)\leq C(\La,
T)\ee^{\frac 14}. \label{z3}\eeq Recall that $\Delta_{\phi}
\dot\phi=n-R(\oo_{\phi})$, by the above estimate and Proposition
\ref{lem2.3} we have \beq|\dot\phi-c(t)|_{C^0}\leq C(\sigma)C(T,
\La)\ee^{\frac 14}, \forall t\in [2T, 6T].\label{z4}\eeq For
simplicity, we  can write $\ee^{\frac 14}$ in the inequalities
(\ref{z3}) and (\ref{z4}) as $\ee$, since we can assume $E_1(0)\leq
\inf E_1+\ee^4$ in the assumption. The lemma is proved.
\end{proof}

\subsection{Estimate of the average of $\pd {\phi}t$}
In this section, we want to control  $c(t)=\frac 1V\int_M\;\dot \phi
\oo^n_{\phi}$. Here we follow the argument in \cite{[chen-tian2]}.
Notice that the argument essentially needs the lower bound of the
$K$-energy, which can be obtained by Theorem \ref{main4} in our
case.  Observe that for any solution $\phi(t)$ of the K\"ahler-Ricci
flow,
$$\pd {\phi}t=\log\frac {\oo_{\phi}^n}{\oo^n}+\phi-h_{\oo},$$
 the function
$\tilde\phi(t)=\phi(t)+Ce^t$ also satisfies the above equation for
any constant $C$. Since $$\pd {\tilde \phi}{t}(0)=\pd
{\phi}{t}(0)+C,$$ we have $\tilde c(0)=c(0)+C.$ Thus we can
normalize  the solution $\phi(t)$ such that the average of
$\dot\phi(0)$ is any given constant.

The proof of the following lemma will be used in section 5 and 6, so
we include a proof here.

\begin{lem}\label{lem2.6}(cf. \cite{[chen-tian2]})Suppose that the $K
$-energy is bounded from below along the K\"ahler-Ricci flow. Then
we can normalize the solution $\phi(t)$ so that
$$c(0)=\frac 1V\int_0^{\infty}\;e^{-t}\int_M\;|\Na \dot\phi|^2\oo_
{\phi}^n\wedge dt<\infty. $$ Then  for all time $t>0$, we have
$$0<c(t),\;\;\int_0^{\infty}\;c(t)dt<E_0(0)-E_0(\infty),$$
where $E_0(\infty)=\lim_{t\ri \infty}E_0(t)$.
\end{lem}
\begin{proof} A simple calculation yields
$$c'(t)=c(t)-\frac 1V\int_M\;|\Na \dot\phi|^2\oo_{\phi}^n.$$
Define
$$\ee(t)=\frac 1V\int_M\;|\Na \dot\phi|^2\oo_{\phi}^n.$$
Since the $K$-energy has a lower bound along the K\"ahler-Ricci
flow, we have
$$\int_0^{\infty}\;\ee(t)dt=\frac 1V\int_0^{\infty}\int_M\;|\Na \dot
\phi|^2\oo_{\phi}^n\wedge dt= E_0(0)-E_0(\infty).$$ Now we normalize
our initial value of $c(t)$ as
\beqs c(0)&=&\int_0^{\infty}\;\ee(t)e^{-t}dt\\
&=&\frac 1V\int_0^{\infty}\;e^{-t}\int_M\;|\Na
\dot\phi|^2\oo_{\phi}^n
\wedge dt\\
&\leq &\frac 1V\int_0^{\infty}\int_M\;|\Na \dot\phi|^2\oo_
{\phi}^n\wedge dt\\
&= &E_0(0)-E_0(\infty). \eeqs  From the equation for $c(t)$, we
have
$$(e^{-t}c(t))'=-\ee(t)e^{-t}.$$
Thus, we have \beqs0<c(t)=\int^{\infty}_t
\;\ee(\tau)e^{-(\tau-t)}d\tau \leq E_0(0)-E_0(\infty) \eeqs and
$$\lim_{t\ri \infty}c(t)=\lim_{t\ri \infty}\int^{\infty}_t \;\ee(\tau)
e^{-(\tau-t)}d\tau=0.$$ Since the $K$-energy is bounded from below,
we have
$$\int_0^{\infty}\;c(t)dt=\frac 1V\int_0^{\infty}\int_M\;|\Na \dot
\phi|^2\oo_{\phi}^n\wedge dt-c(0)\leq E_0(0)-E_0(\infty).$$
\end{proof}

\begin{lem}\label{lem5.9}Suppose that  $E_1$ is bounded from below on
$\cP(M, \oo)$. For any solution $\phi(t)$ of the K\"ahler-Ricci flow
with the initial metric $\oo_0$ satisfying
$$E_1(0)\leq \inf E_1+\ee, $$
after normalization for the K\"ahler potential $\phi(t)$ of the
solution, we have
$$0<c(t),\;\;\int_0^{\infty}c(t)\oo^n_{\phi}\leq \frac {\ee}2.$$
\end{lem}
\begin{proof}By Theorem \ref{main4}, the $K$-energy is bounded from
below, then one can find a sequence of times $t_m\ri \infty$ such
that $$\int_M\; |\Na \dot\phi|^2\oo^n_{\phi}\Big|_{t=t_m}\ri 0.$$ By
Theorem \ref{pali}, we have
$$E_1(t)=2E_0(t)+\frac 1V\int_M\; |\Na\dot\phi|^2\oo_{\phi}^n-C_{\oo}.$$
Then \beqs 2(E_0(0)-E_0(t_m))&=&E_1(0)-E_{1}(t_m)-\frac
1V\int_M\;|\Na\dot\phi|^2\oo_ {\phi}^n\Big|_{t=0}+ \frac
1V\int_M\;|\Na\dot\phi|^2\oo_{\phi}^n\Big|_{t=t_m}\\&\leq&\ee+\frac
1V\int_M\;| \Na\dot\phi|^2\oo_{\phi}^n\Big|_{t=t_m}\\&\ri&\ee.\eeqs
Since the $K$-energy is decreasing along the K\"ahler-Ricci flow, we
have
$$E_0(0)-E_0(\infty)\leq \frac {\ee}2.$$
By the proof of Lemma \ref{lem2.6}, for any solution of the
K\"ahler-Ricci flow we can normalize $\phi(t)$ such that
$$0<c(t),\;\;\int_0^{\infty}c(t)\oo^n_{\phi}\leq E_0(0)-E_0(\infty)
\leq \frac {\ee}2. $$ The lemma is proved.
\end{proof}

\subsection{Estimate of the first eigenvalue of the Laplacian operator}
\subsubsection{Case 1: $M$ has no nonzero holomorphic vector
fields}\label{section4.4.1} In this subsection, we will estimate the
first eigenvalue of the Laplacian when $M$ has no nonzero
holomorphic vector fields. In order to show that the norms of $\phi$
decay exponentially in section 4.5, we need to prove that the first
eigenvalue is strictly greater than $1$.

\begin{theo}\label{lem2.8}Assume that $M$ has no nonzero
holomorphic vector fields. For any $A, B>0$,  there exist $\eta(A,
B, \oo)>0$ such that for any metric
$\oo_{\phi}=\oo+\sqrt{-1}\partial\bar\partial \phi, $ if
$$-\eta \oo_{\phi}\leq Ric(\oo_{\phi})-\oo_{\phi}\leq A\oo_{\phi}\an
|\phi|\leq B,
$$
then the first eigenvalue of the Laplacian $\Delta_{\phi}$ satisfies
$$\la_1>1+\ga(\eta, B, A, \oo),$$
where $\ga>0$ depends only on $\eta, B, A$ and the background
metric $ \oo$.
\end{theo}

The following lemma is taken from \cite{[chen-tian2]}.
\begin{lem}\label{lem2.9}(cf. \cite{[chen-tian2]}) If the K\"ahler
metric $\oo_{\phi}$ satisfies
$$Ric(\oo_{\phi})\geq \al\oo_{\phi},\an |\phi|\leq B$$
for two constants $\al$ and $B,$ then there exists a uniform
constant $C$ depending only on $\al, B$ and $\oo$ such that
$$\inf_M \log \frac {\oo_{\phi}^n}{\oo^n}(x)\geq -4C(\al, B, \La) e^{2
(1+\int_M\; \log \frac {\oo_{\phi}^n}{\oo^n}\oo^n_{\phi})}.$$
\end{lem}

The following crucial lemma is taken from Chen-He \cite{[chenhe]}.
Here we include a proof.
\begin{lem}\label{lem2.10}For any constant $A, B>0$, if
$|Ric(\oo_{\phi})|\leq A$ and $|\phi|\leq B,$ then there is a
constant $C$ depending only on $A, B$ and the background metric
$\oo$ such that $|\phi|_{C^{3, \bb}(M, \oo)}\leq C(A, B, \oo, \bb)$
for any $\bb\in (0, 1)$. In particular, one can find two constants
$C_2(A, B, \oo)$ and $C_3(A, B, \oo)$ such that
$$C_2(A, B,
\oo)\oo\leq \oo_{\phi}\leq C_3(A, B, \oo)\oo.$$

\end{lem}
\begin{proof} We use Yau's estimate on complex Monge-Amp\`ere
equation to obtain the $C^{3, \bb}$ norm of $|\phi|.$ Let $F=\log
\frac {\oo^n_{\phi}}{\oo^n}.$ Then we have \beqs \Delta_{\oo}
F=g^{i\bar j}\partial_i\partial_{\bar j}\log \frac
{\oo^n_{\phi}}{\oo^n} =-g^{i\bar j}R_{i\bar j}(\phi)+R(\oo), \eeqs
where $\Delta_{\oo}$ denotes the Laplacian of $\oo$. On the other
hand, we choose normal coordinates at a point such that $ g_{i\bar
j}=\dd_{ij}$ and $g_{i\bar j}(\phi)=\la_i\dd_{ij},$ then \beqs
g^{i\bar j}R_{i\bar j}({\phi})=\sum_i R_{i\bar i}(\phi) \leq A\sum_i
g_{i\bar i}(\phi) =A(n+\Delta_{\oo} \phi) \eeqs and
$$g^{i\bar j}R_{i\bar j}(\phi)\geq -A(n+\Delta_{\oo}\phi).$$
Hence, we have
\beqn \Delta_{\oo} (F-A\phi)&\leq &R(\oo)+An\label{f1}\\
\Delta_{\oo} (F+A\phi)&\geq &R(\oo)-An. \label{f2}\eeqn Appling the
Green formula, we can bound $F$ from above. In fact, \beqs
F+A\phi&\leq &\frac 1V\int_M\; -G(x, y)\Delta_{\oo}
(F+A\phi)(y)\oo^n
(y)+\frac 1V\int_M\;(F+A\phi) \oo^n\\
&\leq &\frac 1V\int_M\; -G(x, y)(R(\oo)-An)\oo^n(y)+\frac 1V\int_M\;
(F+A\phi) \oo^n\\
&\leq &C(\La, A, B), \eeqs where $\La$ is an upper bound of
$|Rm|_{\oo}$. Notice that in the last inequality we used
$$\frac 1V\int_M\; F\oo^n\leq \log \Big(\frac 1V\int_M\; e^F \oo^n\Big)=0.$$
Hence, $F\leq C(\La, A, B).$ Consider complex Monge-Amp\`ere
equation \beq(\oo+\pbp\phi)^n=e^F\oo^n,\label{f3}\eeq by Yau's
estimate we have \beqs
\Delta_{\phi}(e^{-k\phi}(n+\Delta_{\oo}\phi))&\geq & e^{-k\phi}
(\Delta_{\oo} F-n^2\inf_{i\neq j}R_{i\bar ij\bar j}(\oo))\\
&-&ke^{-k\phi}n(n+\Delta_{\oo}\phi)+(k+\inf_{i\neq j}R_{i\bar ij\bar
j}(\oo))e^{-k\phi+
\frac {-F}{n-1}}(n+\Delta_{\oo}\phi)^{1+\frac 1{n-1}}\\
&\geq &e^{-k\phi}(R(\oo)-An-\Delta_{\oo}\phi-n^2\inf_{i\neq
j}R_{i\bar ij
\bar j}(\oo))\\
&-&ke^{-k\phi}n(n+\Delta_{\oo}\phi)+(k+\inf_{i\neq j}R_{i\bar ij\bar
j} (\oo))e^{-k\phi+ \frac {-F}{n-1}}(n+\Delta_{\oo}\phi)^{1+\frac
1{n-1}}. \eeqs The function $e^{-k\phi}(n+\Delta_{\oo}\phi)$ must
achieve its maximum at some point $p.$ At this point,
$$0\geq -An -\Delta_{\oo}\phi(p)-kn(n+\Delta_{\oo}\phi)+(k-\La)e^{\frac {-F(p)}
{n-1}}(n+\Delta_{\oo}\phi)^{1+\frac 1{n-1}}(p).$$ Notice that we can
bound $\sup F$ by $C(\La, A, B).$ Thus, the above inequality implies
$$n+\Delta_{\oo}\phi \leq C_4(\La, A, B).$$
Since we have an upper bound on $F$, the lower bound of $F$ can be
obtained by Lemma \ref{lem2.9} \beqs\inf F\geq -4C(\La, A, B) \exp
(2+2\int_M\; F \oo_{\phi}^n) =C(\La, A, B).\eeqs On the other hand,
\beqs \inf F\leq\log \frac {\oo^n_{\phi}}{\oo^n}=\log\prod_i
\;(1+\phi_{i\bar i})\leq \log (\prod_i(n+\Delta_{\oo}
\phi)^{n-1}(1+\phi_{i\bar i})). \eeqs Hence, $1+\phi_{i\bar i}\geq
C_5(\La, A, B)>0.$ Thus,
$$C_5(\La, A, B)\leq n+\Delta_{\oo} \phi \leq C_4(\La, A, B).$$
By (\ref{f1}) and (\ref{f2}), we have
$$|\Delta_{\oo} F|\leq C(A, B, \La).$$
By the elliptic estimate, $F\in W^{2, p}(M, \oo)$ for any $p>1.$
Recall that $F$ satisfies the equation (\ref{f3}), we have the
H\"older estimate $\phi\in C^{2, \al}(M, \oo)$  for some $\al\in (0,
1)$ (cf. \cite{[Siu]},\cite{[Tru]}). Let $\psi$ be a local potential
of $\oo$ such that $\oo=\pbp \psi$. Differential the equation
(\ref{f3}), we have
$$\Delta_{\phi}\pd {}{z^i}(\phi+\psi)-\pd {}{z^i}\log \oo^n=\pd {F}{z^i}\in W^{1, p}(M, \oo).$$
Note that the coefficients of $\Delta_{\phi}$ is in $C^{\al}(M,
\oo)$, by the elliptic estimate $\phi\in W^{4, p}(M, \oo)$. Then by
the Sobolev embedding theorem for any $\bb\in (0, 1),$
$$|\phi|_{C^{3, \bb}(M, \oo)}\leq C(A, B, \oo, \bb).$$
The lemma is proved.

\end{proof}

For convenience, we introduce the following definition.
\begin{defi} For any K\"ahler metric $\oo,$ we define
$$W(\oo)=\inf_f \Big\{\int_M \;|f_{\al\bb}|^2\oo^n\;\;\Big|\;f\in W^
{2,2}(M, \oo), \int_M\;f^2\oo^n=1, \int_M\;f\oo^n=0\Big\}.$$
\end{defi}
Assume that $M$ has no nonzero holomorphic vector fields,  then the
following lemma gives a positive lower bound of $W(\oo).$

\begin{lem}\label{lem2.12}Assume that $M$ has no nonzero
holomorphic vector fields. For any constant $A, B>0$, there exists
a positive constant $C_6$ depending on $A, B$ and the background
metric $\oo$, such that for any K\"ahler metric
$\oo_{\phi}=\oo+\pbp\phi,$ if
$$|Ric(\oo_{\phi})|\leq A, \an |\phi|\leq B,$$
then
$$W(\oo_{\phi})\geq C_6>0.$$
\end{lem}
\begin{proof}Suppose not, we can find a sequence of metrics $\oo_m=\oo
+\pbp\phi_m$ and functions $f_m$ satisfying
$$|Ric(\oo_m)|\leq A, \qquad |\phi_m|\leq B,$$
and
$$\int_M\;f_m^2\oo_{m}^n=1, \int_M\;f_m\oo_{m}^n=0, \int_M \;|f_{m,
\al\bb}|_{g_m}^2\oo_{m}^n\ri0.$$ Note that the Sobolev constants
with respect to the metrics $\oo_m$ are uniformly bounded. By Lemma
\ref{lem2.10}, we can assume that $\oo_m$ converges to a K\"ahler
metric $\oo_{\infty}$ in $C^{1, \bb}(M, \oo)$ norm for some $\bb\in
(0, 1).$ Now define a sequence of vector fields \beq
X_m^i=g_m^{i\bar k}\pd {f_{m}}{z^{\bar k}},\qquad X_m=X_m^i\pd
{}{z^i}.\label{f4}\eeq By direct calculation, we have
$$|X_m|^2_{g_m}=|\Na f_m|_{g_m}^2,$$
and \beqs \Big|\pd {X_m}{\bar z}\Big|_{g_m}^2=\sum_{i, j}\Big|\pd
{X_m^i}{z^{\bar j}}\Big|_{g_m}^2=|f_{m, \al\bb}|_{g_m}^2. \eeqs Then
\beq\int_M \;\Big|\pd {X_m}{\bar
z}\Big|_{g_m}^2\oo_{g_m}^n\ri0.\label{la1}\eeq Next we claim that
there exist two positive constants $C_7$ and $C_8$ which depend only
on $A$ and the Poincar\'e constant $\sigma$ such that
\beq0<C_7(\sigma)\leq \int_M \; |X_m|_{g_m}^2 \oo_{g_m}^n\leq
C_8(A).\label{la2}\eeq In fact, since the Poincar\'e constant is
uniformly bounded in our case,
$$\int_M \; |X_m|_{g_m}^2 \oo_{g_m}^n=\int_M \; |\Na f_m|_{g_m}^2 \oo_
{g_m}^n\geq C(\sigma)\int_M\;f_m^2\oo_{g_m}^n=C(\sigma).$$ On the
other hand, since the Ricci curvature has a upper bound, we have
\beqs\int_M \; |\Delta_m f_m|_{g_m}^2 \oo_{g_m}^n&=&\int_M \;|f_{m,
\al \bb}|_{g_m}^2\oo_{g_m}^n+\int_M \;R_{i\bar j}f_{m, \bar i}f_{m,
j}\oo_
{g_m}^n\\
&\leq &\int_M \;|f_{m, \al\bb}|_{g_m}^2\oo_{g_m}^n+A\int_M \;|\Na
f_m|
_{g_m}^2\oo_{g_m}^n\\
&\leq &\int_M \;|f_{m, \al\bb}|_{g_m}^2\oo_{g_m}^n+\frac 12\int_M
\; |
\Delta_m f_m|_{g_m}^2 \oo_{g_m}^n+\frac {A^2}2 \int_M \;f_m^2\oo_{g_m}^n\\
\eeqs Then
$$\int_M \; |\Delta_m f_m|_{g_m}^2 \leq 1+A^2.$$
Therefore, \beqs\int_M \; |X_m|_{g_m}^2 \oo_{g_m}^n&=&\int_M \;
|\Na f_m|_{g_m}
^2 \oo_{g_m}^n\\
&\leq&\frac 12\int_M \; |\Delta_m f_m|_{g_m}^2 + \frac
12\int_M\;f_m^2
\oo_{g_m}^n\\
&\leq &C(A). \eeqs This proves the claim.

Now we have
$$\int_M\;f_m^2\oo_{m}^n=1, \qquad \int_M\;|\Na\bar \Na f_m|_{g_m}^2
\oo_{m}^n\leq C(A),\; \int_M\; |f_{m, \al\bb}|_{g_m}^2 \oo_m^n\ri
0,$$ then $f_m\in W^{2,2}(M, \oo_m).$ Note that the metrics $\oo_m$
are $C^{1, \bb}$ equivalent to $\oo_{\infty}$,  then $f_m\in W^{2,
2}(M, \oo_{\infty}),$ thus we can assume $f_m$ strongly converges to
$f_{\infty}$ in $W^{1, 2}(M, \oo_{\infty})$. By (\ref{f4}) $X_{m}$
strongly converges to $ X_{\infty}$ in $L^2(M, \oo_{\infty}).$ Thus,
by (\ref{la2}), \beq 0<C_7\leq \int_M \; |X_{\infty}|^2
\oo_{\infty}^n\leq C_8.\label{f5}\eeq Next we show that $X_{\infty}$
is holomorphic. In fact, for any vector valued smooth function
$\xi=(\xi^1, \xi^2, \cdots, \xi^n),$ \beqs \Big|\int_M\;
\xi\cdot\bar
\partial X_m\oo_{\infty}^n \Big|^2&=&\Big|\int_M\; \xi^k \frac
{\partial X_m}{\partial \bar z^{ k}}\;\oo_{\infty}^n
\Big|^2\\&\leq &\int_M\; |\xi|^2\oo_{\infty}^n\int_M\; \Big|\pd
{X_m}{\bar z}\Big|^2\oo_{\infty}^n\\&\leq &C \int_M\;
|\xi|^2\oo_{\infty}^n\int_M\; \Big|\pd {X_m}{\bar
z}\Big|_{g_m}^2\oo_{g_m}^n \ri 0. \eeqs On the other hand,
$$\int_M\; \xi\cdot\bar
\partial X_m\oo_{\infty}^n =-\int_M\; \bar \partial \xi\cdot X_m \;
\oo_{\infty}^n\ri -\int_M\; \bar \partial \xi\cdot X_{\infty} \;
\oo_{\infty}^n.$$ Then $X_{\infty}$ is a weak holomorphic vector
field, thus it must be holomorphic. By (\ref{f5}) $X_{\infty}$ is a
{{nonzero}} holomorphic vector field, which contradicts  the
assumption that $M$ has no nonzero holomorphic vector fields. The
lemma is proved.
\end{proof}

\begin{lem}\label{lem2.13} If the K\"ahler metric $\oo_g$ satisfies
$Ric(\oo_g)\geq (1-\eta)\oo_g$ where $0<\eta< \frac {\sqrt{C_6}}2$.
Here $C_6$ is the constant obtained in Lemma \ref{lem2.12}. Then the
first eigenvalue of $\Delta_g$ satisfies $\la_1\geq 1+\ga$, where
$\ga=\frac {\sqrt{C_6}}{2}.$
\end{lem}
\begin{proof} Let $ u$ is any eigenfunction of $\oo_g$ with
eigenvalue $\la_1,$ so $\Delta_g u=-\la_1u.$ Then by direct
calculation, we have \beqs \int_M\; u_{ij}u_{\bar i\bar
j}\;\oo_g^n&=&-\int_M\;u_{ij\bar j}
u_{\bar i}\;\oo_g^n\\
&=& -\int_M\;(u_{j\bar ji}+R_{i\bar k}u_{k})u_{\bar i}\;\oo_g^n\\
&=&\int_M ((\Del_g u)^2-R_{i\bar j}u_{j}u_{\bar i})\;\oo_g^n.
\eeqs This implies \beqs C_6\int_M\; u^2 \oo_g^n&\leq & \int_M\;
((\Del_g u)^2-R_{i\bar
j}u_{j}u_{\bar i})\;\oo_g^n\\
&\leq &\la_1^2 \int_M\; u^2 \oo^n-(1-\eta)\int_M\;|\Na u|^2 \oo_g^n\\
&=&(\la_1^2-(1-\eta)\la_1)\int_M\; u^2 \oo_g^n. \eeqs Thus, we
have $\la_1^2-(1-\eta)\la_1-C_6\geq 0.$ Then,
$$\la_1\geq 1+\frac {\sqrt{C_6}}{2}.$$ \end{proof}

\begin{flushleft}
\begin{proof}[Proof of Theorem \ref{lem2.8}]

The theorem follows directly from the above Lemma \ref{lem2.12} and
\ref{lem2.13}.

\end{proof}
\end{flushleft}

\subsubsection{Case 2: $M$ has nonzero holomorphic vector
fields}\label{section4.4.2} In this subsection, we will consider the
case when $M$ has nonzero holomorphic vector fields. Denote by
$Aut(M)^{\circ}$ the connected component containing the identity of
the holomorphic transformation group of $M$. Let $K$ be a maximal
compact subgroup of $Aut(M)^{\circ}$. Then there is a semidirect
decomposition of $Aut(M)^{\circ}$(cf. \cite{[FM]}),
$$Aut(M)^{\circ}=Aut_r(M)\propto R_u,$$
where $Aut_r(M)\subset Aut(M)^{\circ}$ is a reductive algebraic
subgroup and the complexification of $K$, and $R_u$ is the unipotent
radical of $Aut(M)^{\circ}$.  Let $\eta_r(M, J)$ be the Lie algebra
of $Aut_r(M, J).$

Now we introduce the following definition which is a mild
modification from \cite{[chen1]} and \cite{[PhSt]}.
\begin{defi}\label{prestable} The complex structure $J$ of $M$
is called pre-stable, if no complex structure of the orbit of
diffeomorphism group contains larger (reduced) holomorphic
automorphism group (i.e., $Aut_r(M)$).

\end{defi}

Now we recall the following $C^{k, \al}$ convergence theorem of a
sequence of K\"ahler metrics, which is well-known in literature (cf.
\cite{[PhSt]}, \cite{[Tian4]}).

\begin{theo}\label{conv1} Let $M$ be a compact K\"ahler manifold. Let
$(g(t), J(t))$ be any sequence of metrics $g(t)$ and complex
structures $J(t)$ such that $g(t)$ is K\"ahler with respect to
$J(t)$. Suppose  the following is true:
\begin{enumerate}\item For some integer $k\geq 1$, $|\Na^lRm|_{g(t)}$ is uniformly bounded for
any integer $l (0\leq l< k)$;

\item The injectivity radii $i(M,
g(t))$ are all bounded from below;

\item There exist two
uniform constant $c_1$ and $c_2$ such that $0<c_1\leq \Vol(M,
g(t))\leq c_2$.
\end{enumerate}
Then there exists a subsequence of $t_j$, and a sequence of
diffeomorphism $F_j: M\ri M$ such that the pull-back metrics $\tilde
g(t_j)=F_j^*g(t_j)$  converge in $C^{k, \al}(\forall \,\al\in (0,
1))$ to a  $C^{k, \al}$ metric $g_{\infty}$. The pull-back complex
structure tensors $\tilde J(t_j)=F_j^*J(t_j)$ converge in $C^{k,
\al}$ to an integral complex structure tensor $\tilde J_{\infty}$.
Furthermore, the metric $ g_{\infty}$ is K\"ahler with respect to
the complex structure $\tilde J_{\infty}$.

\end{theo}
\begin{theo}\label{theo4.18}Suppose $M$ is pre-stable.
For any $\La_0, \La_1>0$,   there exists $\eta>0$ depending only on
$\La_0$ and $\La_1$ such that for any metric $\oo\in 2\pi c_1(M),$
if \beq|Ric(\oo)-\oo|\leq \eta,\;\; |Rm(\oo)|\leq \La_0,\;\;|\Na
Rm(\oo)|\leq \La_1, \label{r1}\eeq then for any smooth function $f$
satisfying
$$ \int_M\; f\oo^n=0  {\;\;{ and}\;\;} Re\left(\int_M\; X(f)\oo^n\right)=0, \qquad
\forall X\in \eta(M, J),$$ we have
$$\int_M\; |\Na f|^2\oo^n>(1+\ga(\eta, \La_0, \La_1))\int_M\; |f|^2\oo^n,$$
where $\ga>0$ depends only on $\eta, \La_0$ and $\La_1.$

\end{theo}
\begin{proof}Suppose not, for any positive numbers $\eta_m\ri 0$,
there exists a sequence of K\"ahler metrics $\oo_m\in 2\pi c_1(M)$
such that \beq |Ric(\oo_m)-\oo_{m}|\leq \eta_m,\;\; |Rm(\oo_m)|\leq
\La_0,\;\;\;|\Na_m Rm(\oo_m)|\leq \La_1,\label{r2}\eeq where $Rm_m$
is with respect to the metric $\oo_m$,  and smooth functions $f_m$
satisfying
$$\int_M\; f_m\oo_m^n=0,  \qquad Re\left(\int_M\;
X(f_m)\oo_m^n\right)=0, \qquad \forall X\in \eta(M, J),$$
\beq\int_M\; |\Na_m f_m|^2\oo_m^n<(1+\ga_m)\int_M\;
|f_m|^2\oo_m^n,\label{eq5.25}\eeq where $0<\ga_m\ri 0.$  Without
loss of generality, we may assume that
\[
\int_M\; f_m^2 \omega_m^n = 1, \qquad \forall m \in \NN,
\]
which means
\[
\int_M\;
|\Na_m f_m|^2\oo_m^n\leq  1 + \gamma_m < 2.
\]
Then, $f_m $ will converge in $W^{1,2}$ if $(M, \omega_m)$
converges.  However, according to our stated condition, $(M,
\omega_m, J)$ will converge in $C^{2, \al}(\al\in (0, 1))$ to $(M,
\omega_\infty, J_\infty).\;$ In fact, by (\ref{r2}) the diameters of
$\oo_m$ are uniformly bounded. Note that all the metrics $\oo_m$ are
in the same K\"ahler class, the volume is fixed. Then by (\ref{r2})
again, the injectivity radii are uniformly bounded from below.
Therefore, all the  conditions of Theorem \ref{conv1} are satisfied.

 Note that the complex structure
$J_\infty$ lies in the closure of the orbit of diffeomorphisms,
while $\omega_\infty$ is a K\"ahler-Einstein metric in $(M,
J_\infty)$.   By the standard deformation theorem in complex
structures, we have
\[
   \dim  Aut_r (M, J) \leq \dim Aut_r(M, J_\infty).
\]
By abusing notation, we can write
\[
Aut_r(M, J) \subset Aut_r(M, J_\infty).
\]
By our assumption of pre-stable of $(M, J)$, we have the
inequality the other way around. Thus, we have
\[
\dim Aut_r(M, J) = \dim Aut_r(M, J_\infty),\;\;\;{\rm or}\;\;\;
Aut_r(M, J) = Aut_r(M, J_\infty).
\]
Now, let $f_\infty $ be the $W^{1,2}$ limit of $f_m$, then we have
\[
  1 \leq  |f_\infty|_{W^{1,2}(M,\, \omega_\infty)} \leq 3
  \]
  and
  \[
  \int_M f_\infty \omega_\infty^n = 0, \qquad Re\left(\int_M\; X(f_\infty) \omega_\infty^n\right) =
  0,\qquad \forall X\in \eta(M, J).
  \]
  Thus, $f_\infty$ is a non-trivial function.  Since $\omega_\infty$ is a K\"ahler-Einstein metric, we have
  \[
  \int_M\; \theta_X f_\infty \omega_\infty^n = 0,
  \]
  where \[
  {\cal L}_X \omega_\infty =\pbp\theta_X.
  \]
  This implies that $f_\infty $ is perpendicular to the first eigenspace\footnote{Note that
   $\triangle \theta_X = -\theta_X $ is totally real for
   $X \in Aut_r(M, J_\infty).\;$ Moreover, the first eigenspace consists of
  all such $\theta_X.\;$} of $\triangle_{\omega_\infty}.\;$
  In other words, there is a $\delta > 0$ such that
  \[
     \int_M |\nabla f_\infty|^2 \omega_\infty^n > (1+ \delta) \int_M f_\infty^2 \omega_\infty^n > 1+ \delta.
  \]
 However, this contradicts the following fact:
 \beqs  \int_M \;|\nabla f_\infty|^2 \omega_\infty^n  & \leq & \displaystyle \lim_{m\rightarrow \infty}
  \int_M |\nabla_m f_m|^2 \omega_m^n \\
 & \leq &  \displaystyle \lim_{m\rightarrow \infty} (1+ \ga_m) \int_M f_\infty^2 \omega_\infty^n = 1.
 \eeqs
The lemma is then proved.
\end{proof}

\subsection{Exponential decay in a short time}\label{section4.5}
In this subsection, we will show that the $W^{1,2}$ norm of $\dot
\phi$ decays exponentially in a short time. Here we follow the
argument in \cite{[chen-tian2]} and use the estimate of the first
eigenvalue obtained in the previous subsection.
\begin{lem}\label{lem2.14} Suppose for any time  $t\in [T_1, T_2]$,
we have
$$|Ric-\oo|(t)\leq C_1\ee\an \la_1(t)\geq 1+\gamma>1.$$
Let $$\mu_0(t)=\frac 1V\int_M\;(\dot\phi-c(t))^2\oo_{\phi}^n.$$ If
$\ee$ is small enough, then there exists a constant $\al_0>0$
depending only on $\gamma, \sigma$ and $C_1\ee$ such that

$$\mu_0(t)\leq e^{-\al_0 (t-T_1)}\mu_0(T_1), \qquad\forall t\in [T_1,
T_2].$$

\end{lem}
\begin{proof} By direct calculation, we have
\beqs \frac {d}{dt}\mu_0(t)&=&\frac
2V\int_M\;(\dot\phi-c(t))(\ddot \phi-c(t)')\oo_{\phi}^n+\frac
1V\int_M\;(\dot \phi-c(t))^2\Delta_{\phi}\dot
\phi\oo_{\phi}^n\\&=&-\frac 2V\int_M\;(1+\dot \phi-c(t))|\Na (\dot
\phi-c(t))|^2\oo_{\phi}^n+\frac 2V\int_M\;(\dot
\phi-c(t))^2\oo_{\phi}^n. \eeqs By the assumption, we have for
$t\in [T_1, T_2]$ \beqs \frac {d}{dt}\mu_0(t)&=&-\frac
2V\int_M\;(1+\dot\phi-c(t))|\Na
\dot\phi|^2\oo_{\phi}^n+\frac 2V\int_M\;(\dot\phi-c(t))^2\oo_{\phi}^n\\
&\leq&-\frac 2V\int_M\;(1-C(\sigma)C_1\ee)|\Na
\dot\phi|^2\oo_{\phi}^n
+\frac 2V\int_M\;(\dot \phi-c(t))^2\oo_{\phi}^n\\
&\leq&-\frac 2V\int_M\;(1-C(\sigma)C_1\ee)(1+\gamma)(\dot
\phi-c(t))^2\oo_{\phi}^n+\frac 2V\int_M\;(\dot
\phi-c(t))^2\oo_{\phi}^n\\&=&-\alpha_0\mu_0(t). \eeqs Here
$$\al_0=2(1-C(\sigma)C_1\ee)(1+\gamma)-2>0,$$ if we choose $\ee$ small
enough. Thus, we have
$$\mu_0(t)\leq e^{-\al_0 (t-T_1)}\mu_0(T_1).$$
\end{proof}

\begin{lem}\label{lem2.15}Suppose for any time $t\in [T_1, T_2]$, we
have
$$|Ric-\oo|(t)\leq C_1\ee\an \la_1(t)\geq 1+\gamma>1.$$
Let $$\mu_1(t)=\frac 1V\int_M\;|\Na \dot \phi|^2\oo_{\phi}^n.$$ If
$ \ee$ is small enough, then there exists a constant $\al_1>0$
depending only on $\gamma$ and $C_1\ee$ such that
$$\mu_1(t)\leq e^{-\al_1(t-T_1)}\mu_1(T_1), \qquad\forall t\in [T_1,
T_2].$$
\end{lem}
\begin{proof} Recall that the evolution equation for $|\Na \dot\phi|^2$
is
$$\pd {}t|\Na \dot \phi|^2=\Delta_{\phi} |\Na \dot \phi|^2-|\Na \Na \dot
\phi|^2-|\Na \bar \Na\dot \phi|^2+|\Na \dot \phi|^2.$$ Then for any
time $t\in [T_1, T_2],$ \beqs \frac d{dt}\mu_1(t)&=&\frac
1V\int_M\;(-|\Na \Na \dot \phi|^2- |\Na \bar \Na\dot \phi|^2+|\Na
\dot \phi|^2+|\Na \dot \phi|^2\Delta_{\phi}
\dot\phi)\;\oo_{\phi}^n\\
&\leq &\frac 1V\int_M\;(-\gamma|\Na \dot
\phi|^2+(n-R(\oo_{\phi}))|\Na \dot \phi|^2)
\oo_{\phi}^n\\
&\leq& -(\gamma-C_1\ee)\mu_1(t). \eeqs Thus, we have
$$\mu_1(t)\leq e^{-\al_1(t-T_1)}\mu_1(T_1)$$
where $\al_1=\gamma-C_1\ee>0$ if we choose $\ee$ small.
\end{proof}

\subsection {Estimate of the $C^0$ norm of $\phi(t)$ }\label{section4.6}
In this subsection, we derive some estimates on the $C^0$ norm of
$|\phi|$. Recall that in the previous subsection  we proved that the
$W^{1, 2}$ norm of $|\dot\phi-c(t)|$ decays exponentially. Based on
this result we will use the parabolic Moser iteration to show that
the $C^0$ norm of $|\dot\phi-c(t)|$ also decays exponentially.

\begin{lem}\label{lem2.18} Suppose that $\mu_0(t), \mu_1(t)$ decay
exponentially for $t\in [T_1, T_2]$ as in Lemma \ref{lem2.14} and
\ref{lem2.15}, then we have
$$\Big|\pd {\phi}t-c(t)\Big|_{C^{0}}\leq \frac {C_9(n, \sigma)}
{\tau^{\frac {m}{4}}}\Big(\mu_0(t-\tau)+\frac
1{\al_1^2}\mu_1^2(t-\tau)\Big)^{\frac 12},\;\;\;\forall \;t\in
[T_1+\tau, T_2]$$ where $m=\dim_{\RR}M$ and $\tau<T_2-T_1$.

\end{lem}
\begin{proof} Let $u=\pd {\phi}t-c(t),$ the  evolution equation
for $u$ is
$$\pd ut=\Delta_{\phi} u+u+\mu_1(t),$$
where $\mu_1(t)=\frac 1V\int_M\;|\Na \dot \phi|^2\oo_{\phi}^n.$
Note that in the proof of Lemma \ref{lem2.15}, we derived
$$\pd {}t\mu_1(t)\leq -\al_1 \mu_1(t).$$
Thus, we have
$$\pd {}t(u_++\frac 1{\al_1}\mu_1)\leq \Delta_{\phi} (u_++\frac 1{\al_1}
\mu_1)+(u_++\frac 1{\al_1}\mu_1).$$ where $u_+=\max\{u, 0\}$.
Since $u_++\frac 1{\al_1}\mu_1$ is a nonnegative function, we can
use the parabolic Moser iteration,
$$(u_++\frac 1{\al_1}\mu_1)(t)\leq \frac {C(n, \sigma)}{\tau^{\frac {m
+2}{4}}} \Big(\int_{t-\tau}^t\;\int_M (u_++\frac 1{\al_1}\mu_1)^2(s)
\;\oo_{\phi}^n\wedge ds\Big)^{\frac 12}.$$ Since $\mu_0$ and $\mu_1$
are decreasing,
\beqn& &(u_++\frac 1{\al_1}\mu_1)(t)\nonumber\\
&\leq&\frac {C(n, \sigma)}{\tau^{\frac
{m+2}{4}}}\Big(\int_{t-\tau}^t \;(\mu_0(s)+\frac 1{\al_1^2}
\mu_1^2(s))ds\Big)^{\frac 12}\nonumber\\
&\leq&\frac {C(n, \sigma)}{\tau^{\frac
{m}{4}}}\Big(\mu_0(t-\tau)+\frac
1{\al_1^2}\mu_1^2(t-\tau)\Big)^{\frac 12}. \label{x1}\eeqn On the
other hand, the evolution equation for $-u$ is
$$\pd {}t(-u)=\Delta_{\phi} (-u)+(-u)-\mu_1(t)\leq \Delta_{\phi} (-u)+(-u).$$
Thus,
$$\pd {}t(-u)_+\leq \Delta_{\phi} (-u)_++(-u)_+.$$
By the parabolic Moser iteration, we have \beqn  (-u)_+&\leq
&\frac {C(n, \sigma)}{\tau^{\frac {m+2}{4}}}\Big(\int_{t-\tau}^t\;
\int_M (-u)_+^2\oo_{\phi}^n\wedge ds\Big)^{\frac 12}\nonumber\\
&\leq&\frac {C(n, \sigma)}{\tau^{\frac
{m}{4}}}\mu_0(t-\tau)^{\frac 12}. \label{x2}\eeqn Combining the
two inequalities (\ref{x1})(\ref{x2}), we obtain the estimate
$$\Big|\pd {\phi}t-c(t)\Big|_{C^{0}}\leq \frac {C(n, \sigma)}{\tau^
{\frac {m}{4}}}\Big (\mu_0(t-\tau)+\frac
1{\al_1^2}\mu_1^2(t-\tau)\Big)^{\frac 12}.$$ This proved the
lemma.\end{proof}

\begin{lem}\label{lem2.19}Under the same assumptions as in Lemma \ref{lem2.18},  we have
$$|\phi(t)|\leq |\phi(T_1+\tau)|+\frac {C_{10}(n, \sigma)}{\al \tau^{\frac
{m}{4}}}(\sqrt{\mu_0(T_1)}+\frac 1{\al_1}\mu_1(T_1)) + \tilde
C,\qquad \forall\, t\in [T_1+\tau, T_2].$$ Here $\tilde
C=E_0(0)-E_0(\infty)$ is a constant in Lemma \ref{lem2.6}.

\end{lem}
\begin{proof}
\beqs |\phi(t)|&\leq &|\phi(T_1+\tau)|+\int_{T_1+\tau}^{t}\;
\Big|\pd {\phi(s)}
s-c(s)\Big| ds+\int_{T_1+\tau}^{t}\;c(s) ds\\
&\leq &|\phi(T_1+\tau)|+\frac {C(n, \sigma)}{\tau^{\frac
{m}{4}}}\int_{T_1+\tau} ^{t}\;\Big(\mu_0(s-\tau)+\frac
1{\al_1^2}\mu_1^2(s-\tau)\Big)^{\frac 12}ds +
\tilde C\\
&\leq &|\phi(T_1+\tau)|+\frac {C(n, \sigma)}{\tau^{\frac
{m}{4}}}(\sqrt {\mu_0(T_1)}+\frac
1{\al_1}\mu_1(T_1))\int_{T_1+\tau}^{t}\;e^{-
\al (s- \tau-T_1)}ds +\tilde C\\
&\leq&|\phi(T_1)|+\frac {C(n, \sigma)}{\al \tau^{\frac
{m}{4}}}(\sqrt{\mu_0(T_1)}+\frac 1{\al_1}\mu_1(T_1)) +\tilde C \eeqs
where $\al=\min\{\frac {\al_0}2, \al_1\}$ and $\tilde
C=E_0(0)-E_0(\infty)$ is a constant in Lemma \ref{lem2.6}.
\end{proof}

\subsection{Estimate of the $C^k$ norm of $\phi(t)$ }
In this subsection, we shall obtain  uniform $C^k$ bounds for the
solution $\phi(t)$ of the K\"ahler-Ricci flow
$$\pd {\phi}t=\log\frac {\oo^n_{\phi}}{\oo^n}+\phi-h_{\oo}$$
with respect to any background metric $\oo$. For simplicity,  we
normalize $h_{\oo}$ to satisfy
$$\int_M\; h_{\oo}\;\oo^n=0.$$
The following is the main result in this subsection.
\begin{theo}\label{theoRm}For any positive constants $\La, B>0$ and
small $\eta>0$, there exists a constant $C_{11}$ depending only on
$B, \eta, \La$ and the Sobolev constant $\sigma$ such that if the
background metric $\oo$ satisfies
$$|Rm(\oo)|\leq \La, \qquad |Ric(\oo)-\oo|\leq \eta, $$
and $|\phi(t)|, |\dot\phi(t)|\leq B,$ then
$$|Rm|(t)\leq C_{11}(B, \La, \eta, \sigma).$$
\end{theo}
\begin{proof}Note that $R(\oo)-n=\Delta_{\oo} h_{\oo}$, by the assumption we
have
$$|\Delta_{\oo} h_{\oo}|\leq \eta.$$
Since the Sobolev constant with respect to the metric $\oo$ is
uniformly bounded by a constant $\sigma$, we have
$$|h_{\oo}|_{C^0}\leq C(\sigma)\eta.$$
Now we use Yau's estimate to obtain higher order estimate of
$\phi.$ Define
$$F=\dot\phi-\phi+h_{\oo},$$
then the K\"ahler-Ricci flow can be written as
$$(\oo+\pbp \phi)^n=e^F \oo^n.$$
By Yau's estimate we have \beqs
\Delta_{\phi}(e^{-k\phi}(n+\Delta_{\oo}\phi))&\geq &
e^{-k\phi}(\Delta_{\oo} F-n^2
\inf_{i\neq j}R_{i\bar ij\bar j}(\oo))\\
&-&ke^{-k\phi}n(n+\Delta_{\oo}\phi)+(k+\inf_{i\neq j}R_{i\bar ij\bar
j}(\oo))e^{-k\phi+ \frac {-F}{n-1}}(n+\Delta_{\oo}\phi)^{1+\frac
1{n-1}}. \eeqs Note that \beqs \pd {}t(e^{-k\phi}(n+\Delta_{\oo}
\phi))&=&-k\dot\phi e^{-k
\phi}(n+\Delta_{\oo} \phi)+e^{-k\phi}\Delta_{\oo} \dot\phi\\
&=&-k\dot\phi e^{-k\phi}(n+\Delta_{\oo} \phi)+e^{-k\phi}\Delta_{\oo}
(F+\phi-h_{\oo}).\eeqs Combing the above two inequalities, we have
\beqs (\Delta_{\phi}-\pd {}t)(e^{-k\phi}(n+\Delta_{\oo}\phi))&\geq &
e^{-k\phi}(\Delta_{\oo} h_{\oo}+n-n^2\inf_{i\neq j}R_{i\bar ij\bar
j}(\oo))\\&+&(k\dot\phi-kn-1) e^{-k\phi}(n+\Delta_{\oo}
\phi)\\&+&(k+\inf_{i\neq j}R_{i\bar ij\bar j}(\oo))e^{-k\phi+ \frac
{-F}{n-1}}(n+\Delta_{\oo}\phi)^{1+\frac 1{n-1}}.\eeqs Since $\phi,
\Delta_{\oo} h_{\oo}, |h_{\oo}|, |Rm(\oo)|, \dot\phi$ are bounded,
by the maximum principle we can obtain the following estimate
$$n+\Delta_{\oo} \phi\leq C_{12}(B, \eta, \La, \sigma).$$
By the definition of $F$,
$$\log\frac {\oo^n_{\phi}}{\oo^n}=F\geq -C_{13}(B, \eta, \sigma).$$
On the other hand, we have \beqs\log\frac
{\oo^n_{\phi}}{\oo^n}&=&\log\prod_{i=1}^n(1+\phi_{i\bar i})\leq\log
((n+\Del_{\oo}\phi)^{n-1}(1+\phi_{i\bar i})). \eeqs Thus,
$1+\phi_{i\bar i}\geq e^{-C_{13}}C_{12}^{-\frac 1{n-1}},$i.e.
$C_{14}\oo\leq \oo_{\phi}\leq C_{12}\oo.$ Following Calabi's
computation (cf. \cite{[chen-tian1]},\cite{[Yau]}), we can obtain
the following $C^3$ estimate:
$$|\phi|_{C^3(M, \,\oo)}\leq C_{14}(B, \eta, \La, \sigma).$$
Since the metrics $\oo_{\phi}$ are uniformly equivalent, the flow is
uniform parabolic with $C^1$ coefficients. By the standard parabolic
estimates, the $C^4$ norm of $\phi$ is bounded, and then all the
curvature tensors are also bounded. The theorem is proved.
\end{proof}

\section{Proof of Theorem \ref{main} }
In this section, we shall prove Theorem \ref{main}. This theorem
needs the technical condition that $M$ has no nonzero holomorphic
vector fields, which will be removed in Section \ref{section6}. The
idea is to use the estimate of the first eigenvalue proved in
Section \ref{section4.4.1}.

\begin{theo}Suppose that $M$ has no nonzero holomorphic vector fields
and  $E_1$ is bounded from below in $[\oo].$ For any $\delta, B,
\La>0,$ there exists a small positive constant $\ee(\delta, B,
\La, \oo)>0$ such that for any metric $\oo_0$ in the subspace
$\cA(\delta, B, \La, \ee)$ of K\"ahler metrics
$$\{\oo_{\phi}=\oo+\pbp \phi\;|\; Ric(\oo_{\phi})>-1+\delta,\; |\phi|
\leq B, \;|Rm|(\oo_{\phi})\leq \La, \; E_1(\oo_{\phi})\leq \inf
E_1+\ee \}$$ the K\"ahler-Ricci flow will deform it exponentially
fast to a K\"ahler-Einstein metric in the limit.
\end{theo}

\begin{flushleft}
\begin{proof} Let $\oo_0=\oo+\pbp \phi(0)\in \cA(\dd, B, \La,
\ee)$, where $\ee$ will be determined later. Note that $E_1(0)\leq
\inf E_1 +\ee, $  by  Lemma \ref{lem5.9} we have
$$E_0(0)-E_0(\infty)\leq \frac {\ee}2<1.$$
Here we choose $\ee<2.$ Therefore, we can normalize the
K\"ahler-Ricci flow such that for the normalized solution
$\psi(t)$,\beq 0< c(t), \;\int_0^{\infty} c(t)dt<1,\label{a1}\eeq
where $c(t)=\frac 1V\int_M\;\dot\psi\oo^n_{\psi}.$ Now we give the
details on how to normalize the solution to satisfy (\ref{a1}).
Since $\oo_0=\oo+\pbp \phi(0)\in \cA(\dd, B, \La, \ee)$, by Lemma
{\ref{lem2.10}} we have
$$C_2(\La, B, \oo)\oo\leq \oo_0\leq C_3(\La, B, \oo)\oo.$$
By the equation of K\"ahler-Ricci flow,  we have
$$|\dot\phi|(0)= \Big|\log \frac {\oo_{\phi}^n}{\oo^n}+\phi-h_{\oo}
\Big|_{t=0}\leq C_{16}(\oo, \La, B).$$ Set $\psi(t)=\phi(t)+C_0e^t,$
where
$$C_0=\frac 1V\int_0^{\infty}\;e^{-t}\int_M\; |\Na \dot\phi|^2\;\oo^n_
{\phi}\wedge dt -\frac 1V\int_M\; \dot\phi
\oo_{\phi}^n\Big|_{t=0}.$$ Then (\ref{a1}) holds and
$$|C_0|\leq 1+C_{16},$$
and $$|\psi|(0), \;|\dot\psi|(0)\leq B+1+C_{16}:=B_0.$$ \vskip 10pt

\textbf{STEP 1.}(Estimates for $t\in [2T_1, 6T_1]$).  By Lemma
\ref{lem2.1} there exists a constant $T_1(\delta, \La)$ such that
\beq Ric(t)>-1+\frac {\delta}2,\an |Rm|(t)\leq 2\La, \qquad \forall
t\in [0, 6T_1].\label{5.20}\eeq By Lemma {\ref{lem2.5}} and the
equation (\ref{5.20}), we can choose $\ee$ small enough so that
\beq|Ric-\oo|(t)\leq C_1(T_1, \La)\ee<\frac 12,\qquad \forall t\in
[2T_1, 6T_1],\label{5.21}\eeq and \beq|\dot \psi-c(t)|\leq
C(\sigma)C_1(T_1, \La)\ee<1, \qquad\forall t\in [2T_1,
6T_1].\label{5.22}\eeq Then by the inequality
(\ref{a1})\beq|\dot\psi|(t)\leq 1+|c(t)|\leq 2,\qquad \forall t\in
[2T_1, 6T_1]. \label{5.25}\eeq Note that the equation for $\dot\psi$
is
$$\pd {}t\dot\psi=\Delta_{\psi} \dot\psi+\dot\psi,$$
we have \beq |\dot\psi|(t)\leq |\dot\psi|(0)e^{2T_1}\leq
B_0e^{2T_1},\qquad \forall t\in [0, 2T_1]. \label{5.24} \eeq Thus,
for any $t\in [2T_1, 6T_1]$ we have \beqs|\psi|(t)&\leq
&|\psi|(0)+\int_0^{2T_1}\;
|\dot\psi|ds+\int_{2T_1}^t|\dot\psi|ds\\
&\leq &B_0+2T_1B_0e^{2T_1}+8T_1,\eeqs where the last inequality
used (\ref{5.25}) and (\ref{5.24}). For simplicity, we define
\beqs B_1:&=&B_0+2T_1B_0e^{2T_1}+8T_1+2,\\
B_k:&=&B_{k-1}+2,\qquad 2\leq k\leq 4. \eeqs Then
$$|\dot\psi|(t),\;|\psi|(t)\leq B_1,\qquad \forall t\in [2T_1,
6T_1].$$ By Theorem \ref{theoRm} we have
$$|Rm|(t)\leq C_{11}(B_1, \La_{\oo}, 1),\qquad \forall t\in [2T_1,
6T_1],$$ where $\La_{\oo}$ is an upper bound of curvature tensor
with respect to the metric $\oo,$ and $C_{11}$ is a constant
obtained in Theorem \ref{theoRm}. Set $\La_0=C_{11}(B_4, \La_{\oo},
1)$,  we have
$$|Rm|(t)\leq \La_0, \qquad \forall t\in [2T_1, 6T_1].$$

\vskip 10pt \textbf{STEP 2.}(Estimate for $t\in [2T_1+2T_2,
2T_1+6T_2]$). By STEP 1, we have
$$|Ric-\oo|(2T_1)\leq C_1\ee<\frac 12 \an |Rm|(2T_1)\leq\La_0.$$
By Lemma {\ref{lem2.1}}, there exists a constant $T_2(\frac 12,
\La_0)\in (0, T_1]$ such that
$$|Rm|(t)\leq 2\La_0,\an Ric(t)\geq 0,\qquad \forall t\in [2T_1, 2T_1
+6T_2].$$ Recall that $E_1\leq \inf E_1+\ee ,$ by Lemma \ref{lem2.2}
and Lemma \ref{lem2.5} there exists a constant $C_1'(T_2, \La_0)$
such that
$$|Ric-\oo|(t)\leq C_1'(T_2, \La_0)\ee, \qquad\forall t\in [2T_1
+2T_2, 2T_1+6T_2].$$ Choose $\ee$ small enough so that $C_1'(T_2,
\La_0)\ee<\frac 12.$ Then by Lemma \ref{lem2.5},
$$|\dot\psi-c(t)|_{C^0}\leq C(\sigma)C_1'(T_2, \La_0)\ee,\qquad
\forall t\in [2T_1+2T_2, 2T_1+6T_2].$$ Choose $\ee$ small enough so
that $ C(\sigma)C_1'(T_2, \La_0)\ee<1.$ Thus, we can estimate the
$C^0$ norm of $\psi$ for any $t\in [2T_1+2T_2, 2T_1+6T_2]$
\beqs|\psi(t)|&\leq &|\psi|(2T_1+2T_2)+\Big|\int_{2T_1+2T_2}^t\;
\Big(\pd {\psi}{s}-c(s)\Big)ds\Big|+\Big|\int_0^t\;c(s)ds\Big|\\
&\leq &B_1+ 4T_2C(\sigma)C_1'(T_2, \La_0)\ee+1. \eeqs Choose $\ee$
small enough such that $4T_2C(\sigma)C_1'(T_2, \La_0)\ee<1,$ then
$$|\psi(t)|\leq B_2, \qquad\forall t\in [2T_1+2T_2, 2T_1+6T_2].$$
Since $M$ has no nonzero holomorphic vector fields, applying Theorem
\ref{lem2.8} for the parameters $\eta=C_1'\ee,\;A=1\;, |\phi|\leq
B_4,$ if we choose $\ee$ small enough, there exists a constant
$\gamma(C_1'\ee, B_4, 1, \oo)$ such that the first eigenvalue of the
Laplacian $\Delta_{\psi}$ satisfies
$$\la_1(t)\geq 1+\gamma>1, \qquad \forall t\in [2T_1+2T_2, 2T_1+6T_2].$$

\textbf{STEP 3.} In this step, we want to prove the following
claim:
\begin{claim}For any positive number $S\geq 2T_1+6T_2$, if
$$|Ric-\oo|(t)\leq C_1'(T_2, \La_0)\ee<\frac 12 \an |\psi(t)|
\leq B_3, \qquad\forall t\in [2T_1+2T_2, S],$$ then we can extend
the solution $g(t)$ to $[2T_1+2T_2, S+4T_2]$ such that the above
estimates still hold for $t\in [2T_1+2T_2, S+4T_2]$.
\end{claim}
\begin{proof}By the assumption and Lemma {\ref{lem2.5}}, we have
$$|\dot \psi(t)-c(t)|_{C^0}\leq C(\sigma)C_1'(T_2, \La_0)\ee,
\qquad \forall t\in [2T_1+2T_2, S].$$ Note that in step 2, we know
that $C(\sigma)C_1'(T_2, \La_0)\ee<1.$ Then
$$|\dot\psi|(t)\leq 2,\qquad \forall t\in [2T_1+2T_2, S].$$
Therefore, $|\psi|, |\dot\psi|\leq B_3$. By Theorem \ref{theoRm} and
the definition of $\La_0,$ we have
$$|Rm|(t)\leq \La_0,\qquad \forall t\in [2T_1+2T_2, S].$$
By Lemma {\ref{lem2.1}} and the definition of $T_2$,
$$|Rm|(t)\leq 2\La_0, \;\;Ric(t)\geq 0, \qquad\forall t\in [S-2T_2,
S
+4T_2].$$ Thus, by Lemma \ref{lem2.2} and Lemma \ref{lem2.5} we
have
$$|Ric-\oo|(t)\leq C_1'(T_2, \La_0)\ee,\qquad \forall t\in [S, S
+4T_2],$$ and
$$|\dot\psi-c(t)|_{C^0}\leq C(\sigma)C_1'(T_2, \La_0)\ee, \qquad
\forall t\in [S, S+4T_2].$$ Then we can estimate the $C^0$ norm of
$\psi$ for $t\in [S, S+4T_2],$ \beqs|\psi(t)|&\leq
&|\psi|(S)+\Big|\int_S^{S+4T_2}\; \Big(\pd
{\psi}{s}-c(s)\Big)ds\Big|
+\Big|\int_0^{\infty}\;c(s)ds\Big|\\
&\leq &B_3+4T_2C(\sigma)C_1'(T_2, \La_0)\ee+1 \\
&\leq &B_4. \eeqs Then by Theorem \ref{lem2.8} and the definition of
$\ga$, the first eigenvalue of the Laplacian $\Delta_{\psi}$
$$\la_1(t)\geq 1+\gamma>1, \qquad\forall t\in [2T_1+2T_2, S+4T_2].$$
Note that
$$\mu_0(2T_1+2T_2)=\frac 1V\int_M\;(\dot\psi-c(t))^2\oo_{\psi}^n
\leq(C(\sigma)C_1'\ee)^2 $$ and \beqs\mu_1(2T_1+2T_2)&=&\frac
1V\int_M\;|\Na\dot\psi|^2\oo_
{\psi}^n\\
&=&\frac 1V\int_M\;(\dot\psi-c(t))(R(\oo_{\psi})-n)\oo_{\psi}
^n\\
&\leq&C(\sigma)(C_1'\ee)^2. \eeqs By Lemma \ref{lem2.19}, we can
choose $\ee$ small enough  such that  \beqs|\psi(t)|&\leq&
|\psi(2T_1+3T_2)|+\frac {C(n, \sigma)}{\al T_2^{\frac {m}{4}}}
(\sqrt{\mu_0(2T_1+2T_2)}+\frac 1{\al_1}\mu_1(2T_1+2T_2)) +1\\
&\leq&B_2+\frac {C(n, \sigma)}{\al T_2^{\frac {m}{4}}}(1+\frac
1{\al_1}C_1'\ee)
C(\sigma)C_1'\ee+1\\
&\leq &B_3 \eeqs for $t\in [S, S+4T_2].$ Note that $\ee$ doesn't
depend on $S$ here, so it won't become smaller as $S\ri \infty.$
\end{proof}

\textbf{STEP 4.} By step 3, we know the bisectional curvature is
uniformly bounded and the first eigenvalue $\la_1(t)\geq 1+\eta>1$
uniformly for some positive constant $\eta>0.$ Thus, following the
argument in \cite{[chen-tian2]}, the K\"ahler-Ricci flow converges
to a K\"ahler-Einstein metric exponentially fast. This theorem is
proved.

\end{proof}

\end{flushleft}

\section{Proof of Theorem \ref{main2}}\label{section6}
In this section, we shall use the pre-stable condition to drop the
assumptions that $M$ has no nonzero holomorphic vector fields, and
the dependence of the initial K\"ahler potential. The proof here is
roughly the same as in the previous section, but there are some
differences.

In the STEP 1 of the proof below, we will choose a new background
metric at time $t=2T_1$, so the new K\"ahler potential with respect
to the new background metric at $t=2T_1$ is $0$, and has nice
estimates afterwards.  Notice that all the estimates, particularly
in Theorem \ref{theo4.18} and  \ref{theoRm}, are essentially
independent of the choice of the background metric. Therefore the
choice of $\ee$ will not depend on the initial K\"ahler potential
$\phi(0)$. This is why we can remove the assumption on the initial
K\"ahler potential.

As in Theorem \ref{main},  the key point of the proof is to use the
improved estimate on the first eigenvalue in Section
\ref{section4.4.2} (see Claim \ref{last} below). Since the curvature
tensors are bounded in some time interval, by Shi's estimates the
gradient of curvature tensors are also bounded. Then the assumptions
of Theorem \ref{theoRm} are satisfied, and we can use the estimate
of the first eigenvalue.

Now we state the main result of this section.

\begin{theo}Suppose $M$ is pre-stable, and
$E_1$ is bounded from below in $[\oo]$. For any $\delta, \La>0,$
there exists a small positive constant $\ee(\delta, \La)>0$ such
that for any metric $\oo_0$ in the subspace $\cA(\delta, \La, \oo,
\ee)$ of K\"ahler metrics
$$\{\oo_{\phi}=\oo+\pbp \phi\;|\;
Ric(\oo_{\phi})>-1+\delta, \;  |Rm|(\oo_{\phi})\leq \La,\;
E_1(\oo_{\phi})\leq \inf E_1+\ee \},$$ the K\"ahler-Ricci flow will
deform it exponentially fast to a K\"ahler-Einstein metric in the
limit.
\end{theo}

\begin{flushleft}
\begin{proof} Let $\oo_0\in \cA(\dd, \La, \oo, \ee)$, where $\ee$
will be determined later. By Lemma \ref{lem2.1} there exists a
constant $T_1(\delta, \La)$ such that
$$Ric(t)>-1+\frac {\delta}2 \an |Rm|(t)\leq 2\La, \qquad\forall t\in
[0, 6T_1].$$ By Lemma {\ref{lem2.5}},  we can choose $\ee$ small
enough so that \beq|Ric-\oo|(t)\leq C_1(T_1, \La)\ee<\frac 12,
\qquad\forall t\in [2T_1, 6T_1],\label{8.29}\eeq and \beq|\dot
\phi-c(t)|\leq C(\sigma)C_1(T_1, \La)\ee<1, \qquad\forall t\in
[2T_1, 6T_1].\label{8.30}\eeq

\textbf{STEP 1.}(Choose a new background metric). Let $\un
\oo=\oo+\pbp \phi(2T_1)$ and let $\un \phi(t)$ be the solution to
the following K\"ahler-Ricci flow
$$\left\{\begin{array}{l}
  \pd {\un \phi(t)}t=\log \frac {(\un \oo+\pbp \un\phi)^n}{\un \oo^n}
+\un \phi-h_{\un\oo},\qquad t\geq 2T_1,\\
\un \phi(2T_1)=0.\\
\end{array}
\right.    $$ Here $h_{\un \oo}$ satisfies the following
conditions
$$Ric(\un \oo)-\un \oo=\pbp h_{\un \oo}\an \int_M\; h_{\un \oo} \un
\oo^n=0.$$ Then the metric $\un \oo(t)=\un \oo+\pbp \un\phi(t)$
satisfies
$$\pd {}t\un\oo(t)=-Ric (\un \oo(t))+\un \oo(t)\an \un \oo(2T_1)=\oo+
\pbp \phi(2T_1).$$ By the uniqueness of K\"ahler-Ricci flow, we have
$$\un \oo(t)=\oo+\pbp \phi(t),\qquad \forall t\geq 2T_1.$$
Since the Sobolev constant is bounded and
$$|\Delta_{\un \oo}h_{\un \oo}|=|R(\un \oo)-n|\leq C_1(T_1, \La)\ee,$$
we have $$\Big|\pd {\un \phi}{t}\Big|(2T_1)=|h_{\un \oo}|\leq
C(\sigma)C_1(T_1, \La)\ee.$$ Since $E_1$ is decreasing in our case,
we have
$$E_1(\un \oo)\leq E_1(\oo+\pbp \phi(0))\leq \inf E_1+\ee.$$ By
Lemma \ref{lem5.9}, we have $$E_0(\un \oo)\leq \inf E_0+\frac
\ee2.$$ Thus, by Lemma \ref{lem2.6} we have
$$\frac 1V\int_{2T_1}^{\infty}\; e^{-t}\int_M\; \Big|\Na \pd {\un\phi}
{t}\Big|^2\un\oo(t)^n\wedge dt <\frac {\ee}2<1.$$ Here we choose
$\ee<2.$ Set $\psi(t)=\un \phi(t)+C_0e^{t-2T_1},$ where
$$C_0=\frac 1V\int_{2T_1}^{\infty}\; e^{-t}\int_M\; \Big|\Na \pd {\un
\phi}{t}\Big|^2\un\oo(t)^n\wedge dt-\frac 1V \int_M\; \pd {\un
\phi}t\un\oo(t)^n\Big|_{t=2T_1}.
$$
Then $$|\psi(2T_1)|,\;\; \Big|\pd {\psi}{t}\Big|(2T_1)\leq 2,$$ and
$$0< \un c(t),\;\;\int_{2T_1}^{\infty}\; \un c(t)dt<1,$$
where $\un c(t)=\frac 1V\int_M\; \pd {\psi}{t} \un\oo_{\psi}^n.$
Since
$$|\dot\psi-\un c(t)|=|\dot\phi-c(t)|\leq C(\sigma)C_1(T_1, \La)
\ee,\qquad \forall t\in [2T_1, 6T_1],$$ we have
$$|\dot\psi|(t)\leq 2, \qquad \forall t\in [2T_1, 6T_1],$$
and \beqs |\psi|(t)&\leq& |\psi|(2T_1)+\Big|\int_{2T_1}^{t}\;
(\dot\psi-\un c(t))\Big|+\Big|\int_{2T}^{\infty}\; \un
c(s)ds\Big|\\&\leq&3+4T_1C(\sigma)C_1(T_1, \La)\ee, \qquad \forall
t\in [2T_1, 6T_1].\eeqs Choose $\ee$ small enough such that $4T_1
C(\sigma)C_1(T_1, \La)\ee<1$, and define $B_k=2k+2$. Then
$$|\psi|, |\dot\psi|\leq B_1, \qquad \forall t\in [2T_1, 6T_1].$$
By Theorem \ref{theoRm}, we have
$$|Rm|(t)\leq C_{11}(B_1, 2\La, 1),\qquad \forall t\in [2T_1, 6T_1].$$
Here $C_{11}$ is a constant obtained in Theorem \ref{theoRm}. Let
$\La_0:= C_{11}(B_3, 2\La, 1)$, then
$$|Rm|(t)\leq \La_0,\qquad \forall t\in [2T_1, 6T_1].$$

\vskip 10pt \textbf{STEP 2.}(Estimates for $t\in [2T_1+2T_2,
2T_1+6T_2]$). By step 1, we have
$$|Ric-\oo|(2T_1)\leq C_1(T_1, \La)\ee<\frac 12,\an |Rm|(2T_1)\leq
\La_0.$$ By Lemma {\ref{lem2.1}}, there exists a constant $T_2(\frac
12, \La_0)\in (0, T_1]$ such that  \beq|Rm|(t)\leq 2\La_0,\an
Ric(t)\geq 0,\qquad \forall t\in [2T_1, 2T_1+6T_2].\label{7.28}\eeq
Recall that $E_1\leq \inf E_1+\ee ,$ by Lemma \ref{lem2.2} and Lemma
\ref{lem2.5} there exists a constant $C_1'(T_2, \La_0)$ such that
$$|Ric-\oo|(t)\leq C_1'(T_2, \La_0)\ee, \qquad \forall t\in [2T_1+2T_2,
2T_1+6T_2].$$ Choose $\ee$ small enough so that $C_1'(T_2,
\La_0)\ee<\frac 12$. Then by Lemma \ref{lem2.5},
$$|\dot\psi(t)-\un c(t)|_{C^0}\leq C(\sigma)C_1'(T_2, \La_0)\ee,
\qquad \forall t\in [2T_1+2T_2, 2T_1+6T_2].$$ {{Choose $\ee$ small
such that $C(\sigma)C_1'(T_2, \La_0)\ee<1.$}} Thus, we can estimate
the $C^0$ norm of $\psi$ for any $t\in [2T_1+2T_2, 2T_1+6T_2]$
\beqs| {\psi}(t)|&\leq&|\psi|(2T_1+2T_2)+ \Big|\int_{2T_1+2T_2}^t\;
\Big(\pd { \psi}s-\un
c(s)\Big)ds\Big|+\Big|\int_{2T_1+2T_2}^{\infty}\; \un c(s)ds\Big|\\
&\leq &B_1+4T_2 C(\sigma)C_1'(T_2, \La_0)\ee +1\\&\leq &B_2. \eeqs
Here we choose $\ee$ small enough such that $4T_2 C(\sigma)C_1'(T_2,
\La_0)\ee<1$. Thus, by the definition of $\La_0,$ we have
$$|Rm|(t)\leq \La_0, \qquad\forall t\in [2T_1+2T_2, 2T_1+6T_2].$$

\textbf{STEP 3.} In this step, we want to prove the following
claim:
\begin{claim}\label{last}For any positive number $S\geq 2T_1+6T_2$, if
$$|Ric-\oo|(t)\leq C_1'(T_2, \La_0)\ee<\frac 12, \an |Rm|(t)\leq
\La_0,\qquad \forall t\in [2T_1+2T_2, S],$$ then we can extend the
solution $g(t)$ to $[2T_1+2T_2, S+4T_2]$ such that the above
estimates still hold for $t\in [2T_1+2T_2, S+4T_2]$.
\end{claim}
\begin{proof}  By Lemma
{\ref{lem2.1}} and the definition of $T_2$,
$$|Rm|(t)\leq 2\La_0,\; Ric(t)\geq 0, \qquad\forall t\in [2T_1+2T_2, S
+4T_2].$$ Thus, by Lemma \ref{lem2.2} and Lemma \ref{lem2.5} we have
$$|Ric-\oo|(t)\leq C_1'(T_2, \La_0)\ee,\qquad \forall t\in [S-2T_2, S
+4T_2].$$ Therefore, we have
$$|\dot\psi(t)-\un c(t)|_{C^0}\leq C(\sigma)C_1'(T_2, \La_0)\ee,
\qquad \forall t\in [2T_1+2T_2, S+4T_2].$$ By Theorem \ref{main4}
the $K$-energy is bounded from below, then the Futaki invariant
vanishes. Therefore, we have
$$\int_M\; X(\dot\psi)\un\oo_{\psi}^n=0,\qquad \forall X\in
\eta_r(M, J).$$ By the assumption that $M$ is pre-stable and
Theorem \ref{theo4.18}, if $\ee$ is small enough, there exists a
constant $\ga(C_1'\ee, 2\La_0)$ such that
$$\int_M\; |\Na \dot\psi|\un\oo^n_{\psi}\geq (1+\ga)\int_M\; |
\dot\psi-\un c(t)|^2\un\oo^n_{\psi}.$$ Therefore, Lemma
\ref{lem2.14} still holds, i.e. there exists a constant $\al(\ga,
C_1'\ee, \sigma)>0$ such that for any $t\in [2T_1+2T_2, S+4T_2]$
$$\mu_1(t)\leq \mu_1(2T_1+2T_2)e^{-\al(t-2T_1-2T_2)},$$
and
$$\mu_0(t)\leq \frac {1}{1-C_1'\ee}\mu_1(t)\leq 2\mu_1(2T_1+2T_2)e^{-
\al(t-2T_1-2T_2)}. $$

Then by Lemma \ref{lem2.19}, we can choose $\ee$  small enough such
that  \beqs|\psi(t)|&\leq& |\psi(2T_1+3T_2)|+\frac {C_{10}(n,
\sigma)}{\al T_2^{\frac {m}{4}}}
(\sqrt{\mu_0(2T_1+2T_2)}+\frac 1{\al_1}\mu_1(2T_1+2T_2)) +1\\
&\leq&B_2+\frac {C_{10}(n, \sigma)}{\al T_2^{\frac
{m}{4}}}(1+\frac 1{\al_1}C_1'\ee)
C(\sigma)C_1'\ee+1\\
&\leq &B_3 \eeqs for $t\in [S, S+4T_2].$ By the definition of
$\La_0,$ we have
$$|Rm|(t)\leq \La_0,\qquad \forall t\in [S, S+4T_2].$$
\end{proof}

\textbf{STEP 4.} By step 3, we know the bisectional curvature is
uniformly bounded and the $W^{1, 2}$ norm of $\un {\dot \phi}-\un
c(t)$ decays exponentially.  Thus, following the argument in
\cite{[chen-tian2]}, the K\"ahler-Ricci flow converges to a
K\"ahler-Einstein metric exponentially fast. This theorem is proved.

\end{proof}

\end{flushleft}

\vskip3mm

Xiuxiong Chen,  Department of Mathematics, University of
Wisconsin-Madison, Madison WI 53706, USA; xxchen@math.wisc.edu\\

Haozhao  Li, School of Mathematical Sciences, Peking University,
Beijing, 100871, P.R. China; lihaozhao@gmail.com\\

Bing  Wang, Department of Mathematics, University of
Wisconsin-Madison, Madison WI 53706, USA; bwang@math.wisc.edu\\


\begin{thebibliography}{2}
\bibitem{[Pere]}G. Perelman. The entropy formula for the Ricci flow
and its geometric applications. math.DG/0211159.
\end{thebibliography}
\end{document}